\documentclass[12pt]{article}
\usepackage{amsfonts}

\usepackage[T1]{fontenc}

\usepackage{bm}
\usepackage{enumerate} 
\usepackage{amssymb,amsmath}
\usepackage{caption}
\usepackage{subcaption}
\usepackage{accents} 
\usepackage{dsfont}
\usepackage{stackengine}
\let\Horig\H

\usepackage{tikz}
\usetikzlibrary{automata,topaths}
\usetikzlibrary{shapes}
\usetikzlibrary{plotmarks}
 
\usepackage{float} 
\usepackage{textcomp} 
\usepackage{siunitx}
\usepackage{color,cite} 
\definecolor{lightblue}{rgb}{0,0.2,0.5}
\definecolor{blue1}{rgb}{0,0.1,0.9}
\definecolor{mauve}{rgb}{0.7,0,0.43}
\definecolor{dkgreen}{rgb}{0,0.6,0}

\definecolor{codegreen}{rgb}{0,0.6,0}
\definecolor{codegray}{rgb}{0.5,0.5,0.5}
\definecolor{backcolour}{rgb}{0.97,0.97,0.97}

\usepackage[colorlinks=true, urlcolor=lightblue,linkcolor=lightblue, citecolor=lightblue]{hyperref}
\usepackage{ucs}

\usepackage{accsupp}    

\newcommand{\noncopynumber}[1]{
    \BeginAccSupp{method=escape,ActualText={}}
    #1
    \EndAccSupp{}
}

\usepackage{listings}

\lstdefinelanguage{Sage}[]{Python}
{morekeywords={False,sage,True},sensitive=true}

\usepackage{listings}
\makeatletter
\gdef\lst@numberfirstlinefalse{\global\let\lst@ifnumberfirstline\iffalse}
\lst@AddToHook{Init}{\global\let\lst@ifnumberfirstline\iftrue}
\makeatother
\definecolor{dblackcolor}{rgb}{0.0,0.0,0.0}
\definecolor{dbluecolor}{rgb}{0.01,0.02,0.7}
\definecolor{dgreencolor}{rgb}{0.2,0.4,0.0}
\definecolor{dgraycolor}{rgb}{0.30,0.3,0.30}

\usepackage{graphicx}
\usepackage{flushend,cuted}
\usepackage{bm}
\usepackage{tabularx}

\usepackage{indentfirst}
\usepackage{amssymb}
\usepackage{xparse}
\usepackage{tikz}
\usepackage{mdwlist}
\usepackage{tkz-graph}

\DeclareMathAlphabet{\eufrak}{U}{}{}{} 
\SetMathAlphabet\eufrak{normal}{U}{euf}{m}{n}
\SetMathAlphabet\eufrak{bold}{U}{euf}{b}{n}

\newcommand{\R}{\mathbb{R}}

\newcommand{\E}{\mathbb{E}}
\newcommand{\IP}{\mathbb{P}}

\newcommand{\bone}{{\bf 1}}

\newcommand{\N}{\mathbb{N}}

\makeatletter
\newcommand*\rel@kern[1]{\kern#1\dimexpr\macc@kerna}
\newcommand*\widebar[1]{
  \begingroup
  \def\mathaccent##1##2{
    \rel@kern{0.8}
    \overline{\rel@kern{-0.8}\macc@nucleus\rel@kern{0.2}}
    \rel@kern{-0.2}
  }
  \macc@depth\@ne
  \let\math@bgroup\@empty \let\math@egroup\macc@set@skewchar
  \mathsurround\z@ \frozen@everymath{\mathgroup\macc@group\relax}
  \macc@set@skewchar\relax
  \let\mathaccentV\macc@nested@a
  \macc@nested@a\relax111{#1}
  \endgroup
}
\makeatother

\makeatletter
\DeclareRobustCommand\widecheck[1]{{\mathpalette\@widecheck{#1}}}
\def\@widecheck#1#2{
    \setbox\z@\hbox{\m@th$#1#2$}
    \setbox\tw@\hbox{\m@th$#1
       \widehat{
          \vrule\@width\z@\@height\ht\z@
          \vrule\@height\z@\@width\wd\z@}$}
    \dp\tw@-\ht\z@
    \@tempdima\ht\z@ \advance\@tempdima2\ht\tw@ \divide\@tempdima\thr@@
    \setbox\tw@\hbox{
       \raise\@tempdima\hbox{\scalebox{1}[-1]{\lower\@tempdima\box
\tw@}}}
    {\ooalign{\box\tw@ \cr \box\z@}}}
\makeatother

\newtheorem{prop}{Proposition}[section]
\newtheorem{lemma}[prop]{Lemma}
\newtheorem{definition}[prop]{Definition}
\newtheorem{corollary}[prop]{Corollary}

\newtheorem{theorem}[prop]{Theorem}
\newtheorem{remark}[prop]{Remark}
\newtheorem{example}[prop]{Example}
\newtheorem{counterexample}[prop]{Counterexample}
\newtheorem{assumption}[prop]{Assumption}

\def\({\left(}
\def\){\right)}

\def\[{\left[}
\def\]{\right]}
\def\real{{\mathord{\mathbb R}}}
\def\N{{\mathord{\mathbb N}}}

\newenvironment{Proof}{\removelastskip\par\medskip
\noindent{\em Proof.} \rm}{\penalty-20\null\hfill$\square$\par\medbreak}

\newenvironment{Proofy}{\removelastskip\par\medskip
\noindent{\em Proof} \rm}{\penalty-20\null\hfill$\square$\par\medbreak}

\allowdisplaybreaks

\numberwithin{equation}{section}

\usepackage{fvextra}

\usepackage{slashbox}

\GraphInit[vstyle = Shade]
\usetikzlibrary[intersections,
positioning,
petri,
backgrounds,
fit,
decorations.pathmorphing,
arrows,
arrows.meta,
bending,
calc,
intersections,
through,
backgrounds,
shapes.geometric,
quotes,
matrix,
trees,
shapes.symbols,
graphs,
math,
patterns,
external,
scopes,
matrix,
lindenmayersystems,
shapes.callouts,
shapes.misc,
angles,
shapes.arrows,
shadings]

\usetikzlibrary{arrows,automata,shapes,positioning,decorations.pathmorphing,snakes}

\tikzset{snake it/.style={-stealth,
decoration={snake, 
    amplitude = .4mm,
    segment length = 2mm,
    post length=0.9mm},decorate}}

\usetikzlibrary{matrix,calc}

\newcommand{\convexpath}[2]{
[   
    create hullnodes/.code={
        \global\edef\namelist{#1}
        \foreach [count=\counter] \nodename in \namelist {
            \global\edef\numberofnodes{\counter}
            \node at (\nodename) [draw=none,name=hullnode\counter] {};
        }
        \node at (hullnode\numberofnodes) [name=hullnode0,draw=none] {};
        \pgfmathtruncatemacro\lastnumber{\numberofnodes+1}
        \node at (hullnode1) [name=hullnode\lastnumber,draw=none] {};
    },
    create hullnodes
]
($(hullnode1)!#2!-90:(hullnode0)$)
\foreach [
    evaluate=\currentnode as \previousnode using \currentnode-1,
    evaluate=\currentnode as \nextnode using \currentnode+1
    ] \currentnode in {1,...,\numberofnodes} {
-- ($(hullnode\currentnode)!#2!-90:(hullnode\previousnode)$)
  let \p1 = ($(hullnode\currentnode)!#2!-90:(hullnode\previousnode) - (hullnode\currentnode)$),
    \n1 = {atan2(\y1,\x1)},
    \p2 = ($(hullnode\currentnode)!#2!90:(hullnode\nextnode) - (hullnode\currentnode)$),
    \n2 = {atan2(\y2,\x2)},
    \n{delta} = {-Mod(\n1-\n2,360)}
  in 
    {arc [start angle=\n1, delta angle=\n{delta}, radius=#2]}
}
-- cycle
}

\tikzset{hide labels/.style={every label/.append style={text opacity=0}}}

\begin{document}
\title{
\huge
Normal to Poisson phase transition for subgraph counting in the random-connection model 
} 

\author{
  Qingwei Liu\footnote{School of Mathematics and Statistics, The University of Melbourne, Parkville, VIC 3010, Australia.
Email: \href{mailto:qingwei.liu@unimelb.edu.au}{qingwei.liu@unimelb.edu.au}}
  \qquad
      Nicolas Privault\footnote{Division of Mathematical Sciences, School of Physical and Mathematical Sciences, Nanyang Technological University, 21 Nanyang Link, Singapore 637371, Email: \href{mailto:nprivault@ntu.edu.sg}{nprivault@ntu.edu.sg}
}}

\maketitle

\vspace{-0.5cm}

\begin{abstract} 
  We consider the limiting behavior 
  of the count of subgraphs isomorphic to a graph $G$
  with $m\geq 0$ fixed endpoints (or roots) 
  in the random-connection model, 
  as the intensity $\lambda$ of the underlying Poisson
  point process tends to infinity.   
  When connection probabilities are of order  
  $\lambda^{-\alpha}$ we identify a phase transition phenomenon
  depending on a critical decay rate 
  $\alpha^\ast_m (G)>0$
  such that normal approximation for subgraph counts holds
  when $\alpha \in (0,\alpha^\ast_m (G) )$, and a Poisson limit result
  holds if $\alpha = \alpha^\ast_m (G) $. 
  Our approach relies on cumulant growth rates
  derived by the convex analysis 
  of planar diagrams that enumerate the partitions
  involved in cumulant identities. 
  As a result, by the cumulant method we obtain
  normal approximation results
  with convergence rates in the Kolmogorov
  distance, and a Poisson limit theorem, 
  for subgraph counts. 
\end{abstract}

 
\baselineskip0.7cm

\section{Introduction}
\noindent
 The first instance of a threshold \textcolor{black}{phenomenon} in random graphs
  was observed
  \textcolor{black}{in the Erd{\H o}s-R\'enyi model}
  for the containment of balanced graphs,
  \textcolor{black}{see \cite{ER2}}.
  In \cite[Theorem~1]{bollobas81a},
    a Poisson limit theorem
  for the counting of strictly balanced subgraphs
  was proved at the threshold. 
  Poisson approximation under the total variation distance
  \textcolor{black}{has been proved}
  in \cite{barbour82} via the Stein-Chen method, see also \cite[Chapter~5]{BarbourHolstJanson}.
 Phase transition phenomena for inhomogeneous random graphs have
 been studied
 for the growth rate of the giant component in \cite{BJR07}, 
 and for connectivity thresholds in \cite{DevroyeFraiman14}. 

\medskip

\textcolor{black}{Consider an Erd{\H o}s-R\'enyi random graph on $n$ vertices,
  with independent connection probability $p_n$ such that  $p_n = c / n^{\alpha}$
  for some $c>0$ and $\alpha >0$}. 
Necessary and sufficient conditions for the asymptotic normality of the (normalized)
 \textcolor{black}{count $\widebar{N}_G$ of subgraphs} 
 isomorphic to a fixed graph $G$ were obtained in \cite{rucinski},
 as follows. 
     \begin{itemize}
   \item When 
 \begin{equation}
   \label{er1} 
   0 < \alpha < \min_{H\subset G} \frac{v(H)}{e(H)}, 
\end{equation} 
the normalized subgraph count
$\widebar{N}_G$ converges to a normal random variable,
\textcolor{black}{where $v(H)$ and $e(H)$ respectively denote the
  counts of vertices and edges of any subgraph $H$ of $G$.} 
\item
 When
 $$
 \alpha = \min_{H\subset G} \frac{v(H)}{e(H)},
 $$ 
 the rescaled subgraph count converges 
 to a Poisson distribution provided that the graph $G$ is strictly balanced,
 see Definition~\ref{balgra}-\eqref{i1}, and 
 \begin{equation}
   \label{er2} 
   \lim_{n\to \infty} n p_n^{1/\alpha } = c>0, 
\end{equation} 
 see \cite{bollobas1981b}, \cite{karonski1983a}, and Theorem~3.19 in \cite{JLR}.
\end{itemize}
\noindent 
     \textcolor{black}{More recently,} 
explicit convergence rates in the Wasserstein distance for subgraph counts have been obtained in \cite{BKR}, see also \cite{PS2} and references therein for rates in the Kolmogorov distance using the discrete Malliavin calculus.

\medskip

  Poisson and compound Poisson approximation of subgraph counts
  were \textcolor{black}{obtained together} with convergence rates
  in \cite{coulson16} and \cite{coulson18} for the stochastic block model, 
  which can be viewed as a special case of the graphon-based random-connection
  model.
  Normal approximation results for subgraph counts in graphon-based random graphs were obtained in \textcolor{black}{\cite{KaurRollin21}, \cite{zhangzs}, \cite{bhattacharya23}, and \cite{LiuPrivault25}} via the determination of quantitative bounds and higher order fluctuations.
  
\medskip

 This paper considers the counting of subgraphs isomorphic to a fixed graph in the Poisson random-connection model (RCM) $G_\varphi (\eta )$,
which is a random graph whose vertex set is given by a Poisson point process
$\eta$ with intensity $\Lambda$ on $\R^d$, $d\geq 1$,
 and where every pair of vertices is randomly connected with a location-dependent probability given by a connection function $\varphi:\R^d\times \R^d\to[0,1]$.
 The random-connection model can be regarded as a
 unified framework containing several classical random graph models
 as particular cases. 

\medskip

 In particular, when $\varphi$ is a function of the distance between pairs of points of $\eta$, i.e. $\varphi(x,y):=\phi(\|x-y\|)$ for some measurable function $\phi:\R_+\to [0,1]$, the resulting graph is also known as a {soft} random geometric graph, see \cite{penrose91}, \cite{penrose16}, and \cite{LNS21}. 
 When $\phi$ takes the form $\phi ( u ) = \bone_{\{ u \le \varepsilon \}}$, for some
 $\varepsilon >0$, the random-connection model becomes a random geometric graph, c.f. the monograph \cite{penrosebk}, in which a pair of vertices is connected by an edge if and only if the distance between them is less than the fixed threshold $\varepsilon$. 
For another example, when the underlying point process $\eta$ is a binomial point process and the connection function $\varphi(x,y)\equiv p$ is constant for some $p\in(0,1)$, the resulting graph is the Erd{\H o}s-R\'enyi random graph, c.f. \cite{ER} and \cite{G}. When the connection function $\varphi(x,y)$ is symmetric on $[0,1]\times[0,1]$ it is also called a graphon, and the corresponding random-connection model arises as a limit of dense graph sequences, see \cite{LS06}, \cite{lovasz12}, \cite{zhangzs}, \cite{bhattacharya23}.

\medskip

 In the \textcolor{black}{Poisson} random-connection model, 
 phase transition phenomena  have been observed in
\cite{meester},
where a critical point Poisson intensity parameter
has been identified for the occurrence of percolation.
Poisson approximation results for edge counts and for subgraphs of the same order
were established in \cite{penrose18},
and \textcolor{black}{normal}  approximation result
  have been obtained for subgraph counts 
  in \cite{can2022}.
    In \cite{LiuPrivault}, normal approximation for the subgraph counts in the
   Poisson random-connection model as the intensity of the underlying Poisson point process goes to infinity has been established, together with convergence rates under the Kolmogorov distance, through combinatorial arguments and the cumulant method, see also \cite{thale18}, \cite{jansen}, and \cite{doring}, \cite{schulte-thaele}, \cite{heerten} for moderate deviations. By expressing the cumulants of subgraph counts as sums over non-flat and connected partition diagrams, cumulant growth has been analyzed under different limiting regimes, leading to \textcolor{black}{the asymptotic normality of connected subgraph counts} in the dilute case, and for trees in the sparse case. 

\medskip

In this paper, we derive convergence rates in the Kolmogorov distance for the normal approximation of subgraph counts in the Poisson random-connection model, and we
identify a critical threshold
at which Poisson convergence occurs. To the best of our knowledge,
this is the first time that such a phase transition, from normal to Poisson limit theorems for subgraph counts,
 has been observed in the Poisson random-connection model. 
 In addition, we include the counting of subgraphs containing one
 or more endpoints,
 \textcolor{black}{defined as roots placed as arbitrary deterministic
   locations}. 
 \textcolor{black}{This extends} 
 previous approaches to the counting of rooted subgraphs
 in the Erd{\H o}s-R\'enyi
 and inhomogeneous random graph models, see \cite{RV} and \cite{maugis}. 
  In applications to e.g. wireless networks,
  endpoints can model physical devices
  placed at given fixed locations,
  such as roadside units
  in vehicular networks, see, e.g., 
  \cite{ng2011}, \cite{zhang2012}, and 
  \cite{giles-privault2_published}. 
 
\medskip

 We investigate the asymptotic behavior of subgraph counts in a Poisson random-connection model. By allowing the connection probability to vary as the intensity of the underlying Poisson point process tends to infinity, we identify a phase transition where the limiting distribution shifts between normal and Poisson, whenever the graph $G$ satisfies a balance condition \eqref{mbalanced}. Similar to the Erdős-Rényi graph, while this balance condition \eqref{mbalanced} may not be necessary for asymptotic normality, we believe it is essential for the emergence of a Poisson limit, \textcolor{black}{see Remark~\ref{counex2}.}
 
\medskip

The counting of subgraphs containing one or more fixed endpoints, which covers the counting of rooted subgraphs
\cite{RV}, \cite{maugis} in the case of a single endpoint,
 has been considered in \cite{LiuPri23b} 
 in the random-connection model. 
 However, the analyses of cumulant growths in \cite{LiuPrivault} and \cite{LiuPri23b}, were 
 constrained by the use of partitions of maximal cardinality,
 which \textcolor{black}{resulted in}
 an incomplete characterization of the limiting regimes that
 ensure asymptotic normality. 
 In the present paper, we overcome those restrictions via
 a detailed analysis of the behavior of
 the count $N_G$ of subgraphs isomorphic
 to a given connected graph $G$ in the random-connection model
 with and without fixed endpoints.

\medskip

 In comparison with \cite{LiuPrivault} and \cite{LiuPri23b}, 
 the present paper provides a unified analysis
 of normal approximation for subgraph counts
 in terms of a single threshold and in the presence of endpoints,
 without restriction to specific dilute and sparse regimes,
 see Theorem~\ref{fdhskfs} and Corollary~\ref{c37}. 
  This analysis is performed in terms of subgraph densities,
  and it covers the Poisson convergence regime which is shown to hold
  at the parameter threshold, see Theorem~\ref{poslim}.

\medskip

 For this, we develop new combinatorial tools
 in the random-connection model, based
 on the convex analysis of subgraph plots
 introduced for the Erdős-Rényi model in \cite{JLR},  
 which we use to analyse the partitions 
 involved in the representation of subgraph count cumulants.
 This yields an exhaustive analysis of asymptotic normalized cumulant
 growth on a random-connection model $G_\varphi (\eta \cup \{y_1,\ldots , y_m\})$ including endpoints $y_1,\ldots , y_m$,
 with intensity measure of the form $\lambda\cdot\mu$,
 $\lambda >0$, in 
 Proposition~\ref{themo-1}, 
 where $\mu$ is a diffuse sigma-finite measure on $\R^d$
 and the connection function $\varphi:\R^d\times \R^d\to [0,1]$
  is rescaled as $\varphi_\lambda:=c_\lambda\cdot \varphi$,
  $\lambda >0$. 

  \medskip

 Growth rates for the cumulants of normalized subgraph counts
 are then obtained on the random-connection model 
 $G_\varphi (\eta \cup \{y_1,\ldots , y_m\})$
 in Theorem~\ref{fdhskfs} under 
 Assumptions~\ref{assm5-0}-\eqref{assm5-01-2}-\eqref{assm5}
 and the balance condition~\eqref{mbalanced}. 
 \textcolor{black}{Given two functions $f(\lambda )$ and $g(\lambda )>0$,} 
 we write 
 \begin{itemize}
    \item
      $f(\lambda )=O(g(\lambda ))$,
      or $f(\lambda ) \lesssim g(\lambda )$, if $\limsup_{\lambda \to\infty} f(\lambda ) / g(\lambda ) <\infty$,
    \item
      $f(\lambda )=\Omega(g(\lambda ))$, or 
      $f(\lambda ) \gtrsim g(\lambda )$, if $\liminf_{\lambda \to\infty} f(\lambda ) / g(\lambda )>0$,
    \item $f(\lambda )\asymp g(\lambda )$ if $f(\lambda )=O(g(\lambda ))$ and $f(\lambda )=\Omega(g(\lambda ))$, 
\item $f(\lambda )\sim g(\lambda )$ if $\lim_{\lambda \to \infty} f(\lambda )/g(\lambda ) = 1$, 
    \item $f(x)\ll g(x)$, or $g(x)\gg f(x)$, if $f(x)\geq 0$ and $f(x)/g(x)\to 0$.
\end{itemize}
 with the convention $0/0=0$. 
 \textcolor{black}{For} $G$ a graph with $v(G)=r+m$ vertices including $m$ endpoints, 
 \textcolor{black}{we have the following consequences of
 Theorem~\ref{fdhskfs}.} 
 \begin{itemize}
 \item
 We show in Corollary~\ref{c37} 
 that when 
    $$
    \displaystyle
    1 \gtrsim c_\lambda \gg  \lambda^{-\min ( r/e(G) , 1 / a_m(G))}, 
    $$
 the normalized subgraph count
 $\widebar{N}_G$ converges to a normal random variable 
  as $\lambda$ tends to infinity, 
  \textcolor{black}{where $a_m(G)$ is defined in \eqref{maxdegree}
  and depends on endpoint connectivity}.
  As a consequence,
  when 
 $c_\lambda$ takes the form 
 $c_\lambda \asymp \lambda^{-\alpha}$
  as $\lambda$ tends to infinity for some $\alpha >0$,
  we extend the thresholds \eqref{er1}-\eqref{er2}
 in \cite{JLR} from the Erd{\H o}s-R\'enyi model
 to the random-connection model, 
 \textcolor{black}{by showing that,}
  under the balance condition~\eqref{mbalanced}, 
 normal approximation 
   holds for the normalized subgraph count
   $\widebar{N}_G$ provided that 
 $$
   0 < \alpha < \alpha^\ast_m (G) : = \min \left( \frac{r}{e(G)} ,
   \frac{1}{a_m(G)} \right),
$$ 
 i.e. 
 \begin{equation}
   \label{fjkld243} 
 0 < \alpha < \alpha^\ast_m (G)= \frac{r}{e(G)}
\end{equation} 
 when $m=0$ or $m=1$.
\item
\textcolor{black}{In Theorem~\ref{thmstatu},}  
 we derive convergence rates under the Kolmogorov distance, 
 together with 
 a moderate deviation principle, concentration inequalities
 and a normal approximation result with Cram\'er correction. 
 In particular,
 \textcolor{black}{when $c_{\lambda}\asymp \lambda^{-\alpha}$, we obtain} 
 the Kolmogorov bounds 
 \begin{eqnarray}
 \nonumber 
 \lefteqn{
   \sup_{x\in\R}|\IP_\lambda (\widebar{N}_G \leq x)-\Phi(x)|}
 \\
 \nonumber
 & 
 \leq & 
 \left\{\begin{array}{ll}
 \displaystyle         \frac{C}{\displaystyle\lambda^{(1 - \alpha a_m(G))/(4r-2)}} & \ \mbox{if }\
          0 < \alpha \leq \displaystyle \frac{r-1}{e(G)-a_m(G)}
          ,
          \medskip
      \\
\displaystyle
          \frac{C}{\displaystyle\lambda^{(r - \alpha e(G))/(4r-2)}} & \ \mbox{if }\
          \displaystyle \frac{r-1}{e(G)-a_m(G)}\leq \alpha < \alpha^\ast_m (G), 
    \end{array}
    \right.
  \end{eqnarray}
 where $\Phi$ is the cumulative distribution of the
 standard normal distribution and $C>0$ is a constant depending only on $r \geq 2$.
 When $G$ has no endpoints $(m=0)$
 and is strongly balanced, we have 
  \begin{eqnarray}
\nonumber
      \sup_{x\in\R}|\IP_\lambda (\widebar{N}_G \leq x)-\Phi(x)|
     \leq
     \left\{\begin{array}{ll}
     \displaystyle
     \frac{C}{\displaystyle\lambda^{1/(4v(G)-2)}} & \ \mbox{if }\
     \displaystyle 0 < \alpha \leq \frac{v(G)-1}{e(G)}
     ,
     \medskip
      \\
      \displaystyle
      \frac{C}{\displaystyle\lambda^{ ( v(G) - \alpha e(G))/(4v(G)-2)}} & \ \mbox{if }\
      \displaystyle
      \frac{v(G)-1}{e(G)} \leq \alpha <
      \alpha^\ast_0 (G) = \frac{v(G)}{e(G)}, 
    \end{array}
    \right.
\end{eqnarray}
  which extends Corollary~7.1 of \cite{LiuPrivault}
  beyond the dilute regime considered there, 
  and also Corollary~7.2 therein
  without restriction to trees. 

  \medskip
\item
  Under the condition
  $a_m(G) r \leq e(G)$,
  Poisson convergence 
 holds for $N_G$
 by Theorem~\ref{poslim}
 in the boundary case
  $$
  \alpha = \alpha^\ast_m (G)= \frac{r}{e(G)}, 
  $$
\item
  Finally, by Theorem~\ref{pt:contain}-\eqref{part-a}, 
  $N_G$
  converges to zero in probability
  if $a_m(G) r \leq e(G)$ and 
  $$
  \alpha > \alpha^\ast_m (G) = \frac{r}{e(G)}. 
  $$
\end{itemize} 
 \begin{remark}
   \label{jkld14}
    \begin{enumerate}[a)] 
     \item 
         In the case of rooted
 subgraph counting, i.e. 
 when $m=1$ with a single endpoint,
 Condition~\eqref{fjkld243} is consistent
 with the Property~(P) page~261 of 
 \cite{RV}
 in the Erd{\H o}s-R\'enyi model,
 and with the 
 asymptotic normality condition 
 in Theorem~1 of \cite{maugis}
 in inhomogeneous random graphs.
\item  In the absence of endpoints $(m=0)$,
  Condition~\eqref{mbalanced}
  means that $G$ should be strongly balanced, and
 \eqref{fjkld243} reads 
 $$
 0 < \alpha < \alpha^\ast_0 (G)= \min_{H\subset G} \frac{v(H)}{e(H)}
 = \frac{v(G)}{e(G)}, 
$$ 
 which coincides with \eqref{er1}, 
 see Definition~\ref{balgra}-\eqref{i1}. 
\item
 We note that in the random-connection model,
 our results require a strong balance condition of the form
 \eqref{mbalanced}, 
 which is not needed in the
 Erd{\H o}s-R\'enyi model,
 see \cite{rucinski},
 \cite{bollobas1981b},
 and is stronger than
 strict balance,
 see Definition~\ref{balgra}-\eqref{i1}. 
\end{enumerate}
\end{remark} 
\noindent 
 Our approach relies on partition diagrams
 introduced in Section~\ref{sec:diag} and later used 
 in Section~\ref{plot1} to arrange the partitions involved
 in cumulant expressions into a planar representation.
 This planar representation method has been introduced in
 \cite{luczak} to study the behaviour of variance of
 subgraph counts in the Erd{\H o}s-R\'enyi model,
 and is extended here to derive 
 cumulant growth rates of all orders
 in the random-connection model.
 In this paper, it is used in Section~\ref{leading}
 to identify the partition diagrams
 that play a leading role in cumulant expressions, 
 and yields cumulant growth rates in Proposition~\ref{themo-1}. 
 
 \medskip 
 
 We proceed as follows.
  After recalling necessary preliminaries on
 the random-connection model and balanced graphs 
 in Section~\ref{prelim},
 we present our main results in Section~\ref{mainres}. 
 Sections~\ref{sec:diag} and \ref{plot1}  
 focus the planar diagram representation of cumulants,
 and Section~\ref{leading} identifies
 the leading diagrams
 appearing in cumulant expressions.
 Finally, growth rates for cumulants
 are derived in Section~\ref{cumulantrates}. 
 SageMath and R codes used for the computation of
convex hulls and for partition counting 
are listed in Appendices~\ref{fjkldsf-11}
and \ref{fjkldsf-22}.

\section{Preliminaries and notation}
\label{prelim}
\subsubsection*{Random-connection model}
 \noindent
 Given $\lambda > 0$ and $\mu$ a diffuse sigma-finite
 measure on $\real^d$, $d\geq 1$, 
 we consider a Poisson point process $\eta$ on $\R^d$ 
 with intensity measure of the form
 $\lambda\cdot\mu$, \textcolor{black}{which} can be almost surely written as 
\begin{equation}
   \nonumber 
 \eta=\sum_{i=1}^\tau\delta_{X_i} 
\end{equation}
under the probability measure
$\IP_\lambda$, see \cite[Corollary~6.5]{LastPenrose17}, 
 where $\tau$ is a $\N\cup\{\infty\}-$valued random variable, $\delta_x$ denotes the Dirac measure at $x\in\R^d$, and $X_1,X_2,\dots$ are random elements in $\R^d$. For fixed $m\geq 0$, and $y_1,\dots,y_m\in\R^d$, 
 we consider the point process $\eta \cup \{y_1,\ldots , y_m\}$ on $\R^d$
 defined by the union of $\eta$ and $y_1,\ldots , y_m$. 

 \medskip

 Given $\varphi:\R^d\times \R^d\to [0,1]$ a symmetric measurable function
 and $c_\lambda\in(0,1)$, 
 we also let $\varphi_\lambda:=c_\lambda\cdot \varphi$,
 $\lambda >0$,
 denote the connection function of \textcolor{black}{the} random-connection model $G_\varphi (\eta \cup \{y_1,\ldots , y_m\})$. 
 The random-connection model 
 is a random graph
 denoted by $G_\varphi (\eta \cup \{y_1,\ldots , y_m\})$,
 with vertex set $\eta \cup \{y_1,\ldots , y_m\}$,
 such that any two distinct vertices $x,y\in\eta \cup \{y_1,\ldots , y_m\}$ are independently connected by an edge with probability $\varphi_\lambda(x,y)$.
\subsubsection*{Balanced graphs}
\noindent
In what follows, for any two graphs $G_1,G_2$, we write $G_1\simeq G_2$ when $G_1$ is isomorphic to $G_2$. 
We also let $v(G):=|V_G| \geq 2$ and $e(G):=|E_G|$ be the number of vertices and the number of edges of any graph $G$.
\begin{definition}{\cite{luczak},  \cite[pages~64-65]{JLR}} 
\label{balgra}
  \begin{enumerate}[1)]
    \item
      \label{i1}
      A graph $G$ is balanced if 
    \begin{equation}
 \label{fjl3kjl} 
\frac{e(H)}{v(H)}\leq \frac{e(G)}{v(G)},
\qquad H\subset G, 
  \end{equation}
 and strictly balanced if \eqref{fjl3kjl} 
 holds as a strict inequality for all $H\subsetneq G$. 
\item 
 A graph $G$ is strongly balanced 
 if 
\begin{equation}\label{strongbaldef}
  \frac{e(H)}{v(H)-1}\leq \frac{e(G)}{v(G)-1},
  \qquad
  H\subset G, 
\end{equation}
and strictly strongly balanced if
\eqref{strongbaldef} holds as a strict inequality for
 all $H\subsetneq G$. 
 \item 
 A graph $G$ is $K_2$-balanced
 if 
\begin{equation}\label{k2baldef}
  \frac{e(H)-1}{v(H)-2}\leq \frac{e(G)-1}{v(G)-2},
  \qquad
  H\subset G,~ v(H)\geq 3,
\end{equation}
and strictly $K_2$-balanced if
\eqref{k2baldef} holds as a strict inequality for
 all $H\subsetneq G$. 
\end{enumerate}
\end{definition}
\begin{remark}\label{graph-example}
  From \cite{luczak} we have the following statements. 
  \begin{enumerate}[i)]
\item
  Cycles and complete graphs are strictly $K_2$-balanced.
\item
  Trees are $K_2$-balanced, but not strictly $K_2$-balanced.
  \item
  $K_2$-balanced graph are strongly balanced,
  except for the unions of disjoint edges, \textcolor{black}{also called matchings}.
\item
 Strongly balanced graphs are strictly balanced. 
  \end{enumerate}
  \end{remark} 

\subsubsection*{Graphs with endpoints}
\noindent
 Throughout this paper, we consider a connected graph $G$
 satisfying the following conditions. 
\begin{assumption}
\label{a0}
 Given $r\geq 2$ and $m\geq 0$, we 
consider a connected graph $G=(V_G,E_G)$ with edge set $E_G$ and
vertex set $V_G=\{1, \ldots , r+m\}$, such that
\begin{enumerate}[i)]
\item the subgraph 
  induced by $G$ on $\{1, \ldots ,r\}$ is connected, and 
\item \label{ii}
 the endpoint vertices $r+1, \ldots ,r+m$ are not adjacent to each other in $G$, 
\end{enumerate}
where Condition~\eqref{ii} is void and
$V_G=\{1, \ldots ,r\}$ in case $m=0$.
\end{assumption}
\begin{figure}[H]
 \captionsetup[subfigure]{font=footnotesize}
\centering
\subcaptionbox{No endpoint, $m=0$.}[.32\textwidth]{
  \begin{tikzpicture}
    \coordinate (A) at (0,0);  
  \coordinate (B) at (2,0);  
  \coordinate (C) at (1,2);  

    \draw[thick] (A) -- (B) -- (C) -- cycle;

  \fill[blue] (A) circle (2pt);
  \fill[blue] (B) circle (2pt);
  \fill[blue] (C) circle (2pt);

\end{tikzpicture}
}
 \captionsetup[subfigure]{font=footnotesize}
\centering
\subcaptionbox{One endpoint (rooted graph).}[.32\textwidth]{
 \begin{tikzpicture}
    \coordinate (A) at (0,0);  
  \coordinate (B) at (2,0);  
  \coordinate (C) at (1,2);  

    \coordinate (D) at (-2,0);

  \draw[thick] (A) -- (B) -- (C) -- cycle;

    \draw[black,thick] (D) -- (A);

  \fill[blue] (A) circle (2pt);
  \fill[blue] (B) circle (2pt);
  \fill[blue] (C) circle (2pt);

    \fill[red] (D) circle (2pt);

\end{tikzpicture}
}
 \captionsetup[subfigure]{font=footnotesize}
\centering
\subcaptionbox{Two endpoints, $m=2$.}[.32\textwidth]{
 \begin{tikzpicture}
  
  \coordinate (A) at (0,0);  
  \coordinate (B) at (2,0);  
  \coordinate (C) at (1,2);

  \coordinate (D) at (-2,0); 
  \coordinate (E) at (-1,2);

  \draw[thick] (A) -- (B) -- (C) -- cycle;

  \draw[black,thick] (D) -- (A);
  \draw[black,thick] (E) -- (C);

  \fill[blue] (A) circle (2pt);
  \fill[blue] (B) circle (2pt);
  \fill[blue] (C) circle (2pt);

    \fill[red] (D) circle (2pt);
  \fill[red] (E) circle (2pt);
\end{tikzpicture}
}
\caption{Examples of triangles with endpoints, $r=3$.} 
\end{figure}
In the sequel, we
denote $[n]:=\{1,\dots,n\}$ for any $n\geq 1$,
and write $V_G=[r+m]$.
    \begin{definition}
      We let 
\begin{equation}\label{maxdegree}
      a_m(G):=\max_{i\in[r]} |A_i|
\end{equation}
denote the maximum number of endpoint connections to any vertex in $[r]$,
where 
    \begin{equation}
      A_i:=\{ j \in \{r+1 , \ldots , r+m \} \ : \ \{i,j\}\in E_G\},
    \end{equation}
    is the neighborhood of vertex
      $i\in[r]$ within the set
    $\{r+1 , \ldots , r+m \}$ of endpoints.
    \end{definition}
    We note that 
$$a_m(G)\leq m, \quad a_0(G)=0, \quad \mbox{and} \quad a_1(G)=1.
$$ 
    In what follows, our main results will
    hold under the balance condition 
 \begin{equation}
   \label{mbalanced}
    \frac{e(H)-a_m(G)}{v(H)-m-1}\leq \frac{e(G)-a_m(G)}{r-1},\qquad H\subset G,~v(H) \geq m+2\textcolor{black}{.} 
 \end{equation}
 We also note the following points. 
\begin{remark}
  \label{jlkdfed23} 
\begin{enumerate}[a)] 
\item
 When $r=2$ and $m\geq 0$, Condition~\eqref{mbalanced} is satisfied by all
 connected graphs.
\item When $r\geq 3$ and $m\geq 0$, Condition~\eqref{mbalanced}
 is satisfied by any tree $G$, if \textcolor{black}{$m=a_m(G)$.}
 Indeed, when $G$ is a tree and $H$ is a subgraph of $G$, we have
$$
 \frac{e(G)-a_m(G)}{r-1}=
 \frac{r-1+m-a_m(G)}{r-1} =
 1 + \frac{m-a_m(G)}{r-1}
$$ 
 and 
$$
 \frac{e(H)-a_m(G)}{v(H)-m-1}\textcolor{black}{\le}
 \frac{v(H)-1-a_m(G)}{v(H)-m-1}=
 1 + \frac{m-a_m(G)}{v(H)-m-1}, 
$$ 
 hence \eqref{mbalanced} is satisfied if $m = a_m(G)$.
\item When $r\geq 2$ and $m=0$, \eqref{mbalanced}
 is the strong balance condition \eqref{strongbaldef}  
 since
 $a_0(G)=0$. 
\item When $r\geq 2$ and $m=1$, \eqref{mbalanced}
    is   the $K_2$-balance condition \eqref{k2baldef}
    since $a_1(G)=1$.
 Note that cycles, trees and complete graphs
 are $K_2$-balanced by Remark~\ref{graph-example}. 
\end{enumerate}
\end{remark}

\noindent
 Table~\ref{table134}
 presents the counts of \textcolor{black}{(isomorphic)} trees ($\eufrak{t}$) {\em vs.}  
 graphs ($\eufrak{g}$) satisfying Condition~\eqref{mbalanced}
 within 
 those ($\eufrak{a}$) satisfying Assumption~\ref{a0}
 in the format 
 $\eufrak{t}/\eufrak{g}/\eufrak{a}$, 
 using the R code presented in Appendix~\ref{fjkldsf-22}
 for different values of $r\geq 2$ and $m\geq 0$,
 with $m+r \leq 8$.
 We check that when $m\geq 1$ the number of trees satisfying
 Condition~\eqref{mbalanced} depends only \textcolor{black}{on} $r\geq 3$, and 
 that only trees can satisfy Condition~\eqref{mbalanced} 
 as the number $m$ of endpoints becomes large.
 The first row ($m=0$) refers to strongly balanced graphs.
 
\begin{table}[H]
\begin{center}
\begin{tabular}{|c|c|c|c|cc}
\hline
  \backslashbox{$m$}{$r$} & 2 & 3 & 4 & \multicolumn{1}{c|} 5 & \multicolumn{1}{c|} 6 \\
  \hline
  0 & 1/1/1 &  1/2/2 &  2/5/6 & \multicolumn{1}{c|}{3/14/21} & \multicolumn{1}{c|} {6/53/112} \\
  \hline
  1 &  1/2/2 & 2/6/8 &  4/20/44 & \multicolumn{1}{c|}{9/106/333} & \multicolumn{1}{c|} {20/893/3771} \\
  \hline
  2 &  2/4/4 & 2/6/27 &  4/26/274 & \multicolumn{1}{c|}{9/176/4071} & \multicolumn{1}{c|}{20/2273/94584} \\
  \hline
  3 &  2/6/6 & 2/2/73 & 4/7/1346 & \multicolumn{1}{c|}{9/27/39159} 
  \\
  \cline{1-5}
  4 & 3/9/9 & 2/2/171 & \multicolumn{1}{c|}{4/4/5620}   
  \\
  \cline{1-4}
  5 & 3/12/12 & 2/2/359 
  \\
  \cline{1-3}
  6 & 4/16/16 
  \\
  \cline{1-2}
\end{tabular}
\caption{
  Counts $\eufrak{t}/\eufrak{g}/\eufrak{a}$
  of graphs $G$ satisfying Condition~\eqref{mbalanced}
 {\em vs}. Assumption~\ref{a0}. } 
\label{table134}
\end{center} 
\end{table}

 \vspace{-1.2cm}
\section{Main results} 
\label{mainres} 
\noindent
Let $y_1,\dots,y_m\in\R^d$ be fixed endpoints, or terminal nodes, where $m\geq 0$
and by convention we set $\{y_1,\dots,y_m\} = \emptyset $ when $m=0$.
 In what follows, we consider the count $N_G$ of subgraphs
 isomorphic to a given connected graph $G$ in the random-connection model $G_\varphi (\eta \cup \{y_1,\ldots , y_m\})$ 
 which includes the fixed endpoints $y_1,\ldots , y_m$ as vertices.
  \begin{definition}
   \label{jkldaa1}
    Let $N_G$ denote the count of \textcolor{black}{(labelled)} subgraphs $H\subset G_\varphi (\eta \cup \{y_1,\ldots , y_m\})$ such that there exists a bijection $\psi:[r+m]\to V_H$ satisfying \textcolor{black}{
   $$
   \{i,j\}\in E_G\quad \mbox{iff} \quad \{\psi(i),\psi(j)\}\in E_H$$}
   for $1\le i\ne j\le r+m$, and $\psi( l )=y_l$, $l=1,\dots,m$.
 \end{definition}
     We note that the subgraph count $N_G$ can be represented as 
\begin{equation}\nonumber 
  N_G:=\sum_{(\alpha_1,\dots,\alpha_r)\in[\tau]^r_{\ne}}
  \left(
  \prod_{\substack{\{i,j\}\in E_G\\i\in[r],~j\in[r]}}\bone_{\{X_{\alpha_i}\leftrightarrow X_{\alpha_j}\}}
  \right)
  \left(
  \prod_{\substack{\{r+j,i\}\in E_G\\i\in[r],~j\in[m]}}\bone_{\{y_j\leftrightarrow X_{\alpha_i}\}}\right) 
   \end{equation}
 where $x\leftrightarrow y$ indicates that
 $x,y\in\eta \cup \{y_1,\ldots , y_m\}$ are connected by an edge in the RCM, and $$[\tau]^r_{\ne}:=\{(i_1,\dots,i_r)\in[\tau]^r: i_k\ne i_j~\mathrm{if}~k\ne j\}.
 $$
 When $m=0$, $N_G$ is the count
 of graphs isomorphic to the graph $G$
 in the Poisson random-connection model $G_\varphi (\eta )$.
\begin{assumption}\label{assm5-0}
   \begin{enumerate}[i)]
   \item
     \label{assm5-01-2}
     When $m=0$, 
     $\mu$ is a finite diffuse measure on \textcolor{black}{$\real^d$}, 
     and $\varphi:\R^d\times \R^d\to [0,1]$ is a symmetric measurable function. 
\item
    \label{assm5-1-2}
    When $m\geq 1$,
    $\mu$ is a
    diffuse sigma-finite measure on $\R^d$, 
   and $\varphi:\R^d\times \R^d\to [0,1]$ is a symmetric measurable function  
 that satisfies 
 \begin{equation}
\label{integ1-0}
 \sup_{x\in\R^d}\int_{\R^d}\varphi(x,y) \ \! \mu(\mathrm{d}y)<\infty. 
\end{equation}
\item
   \label{assm5}
    When $m\geq 1$,
    in addition \textcolor{black}{to} \eqref{assm5-1-2},
    $\mu$ is the Lebesgue measure on $\real^d$ and
    $\varphi:\R^d\times \R^d\to [0,1]$ is a symmetric measurable function
which is
    translation invariant, i.e. 
    $$\varphi(x,y)=\varphi(0,y-x), 
    \quad
    x,y\in\R^d. 
    $$
\end{enumerate}
\end{assumption}
 Under Assumption~\ref{assm5-0}-\eqref{assm5},
 the integrability condition \eqref{integ1-0} 
 reads 
    \begin{equation}
      \int_{\R^d}\varphi(0,x) \ \! \mathrm{d}x<\infty. 
    \end{equation}
\subsubsection*{Normal approximation} 
\noindent
\textcolor{black}{
 Recall that by e.g. Theorem~1 in \cite{Janson1988},
 any sequence $(X_n)_{n\geq 1}$ of real-valued
 such that  
  $$
  \lim_{n\to \infty} \kappa_m (X_n) = 0,
  \quad
  \mbox{for all }
  m \geq m_0,
  $$
  for some $m_0 \geq 3$
  converges in distribution to the Gaussian
  distribution ${\cal N}(\mu , \sigma^2)$,
  provided that the limits
  $$
  \mu : = \lim_{n\to \infty} \kappa_1 (X_n)
  \quad
  \mbox{and}
  \quad
  \sigma^2 : = \lim_{n\to \infty} \kappa_2 (X_n) 
  $$
  exist. 
  Theorem~\ref{fdhskfs}, which 
  is a consequence of Propositions~\ref{pp:mbal} and \ref{thmstatu0},
  provides sufficient conditions for the asymptotic
  vanishing of higher order cumulants in \eqref{cumboun1}. 
  This will further enable us to apply the method of cumulant for normal
  approximation in Theorem~\ref{thmstatu}.} 
\begin{theorem}
  \label{fdhskfs}
        Let $G$ be a connected graph with $V_G=[r+m]$ for $r\geq 2$ and $m\geq 0$,
        and suppose that
        Assumptions~\ref{assm5-0}-\eqref{assm5-01-2}-\eqref{assm5} are satisfied
        and the balance condition~\eqref{mbalanced} holds. 
  Then, the cumulant $\kappa_n(\widebar{N}_G)$ of order $n\geq 1$ of
 the normalized subgraph count 
  $$
 \widebar{N}_G:=\frac{N_G-\kappa_1(N_G)}{\sqrt{\kappa_2(N_G)}}$$
 satisfies the cumulant bound 
  \begin{equation}\label{cumboun1}
    |\kappa_n(\widebar{N}_G)|\leq \frac{n!^r }{\Delta_\lambda^{n-2}}, 
  \end{equation}
  where 
    \begin{eqnarray}
 \label{d1} 
    \Delta_\lambda\asymp\left\{\begin{array}{ll}
          \displaystyle
          \lambda^{r/2}c_\lambda^{e(G)/2} & \ \mbox{if }\
           c_\lambda
      \lesssim \lambda^{-(r-1)/(e(G)-a_m(G))},
          \medskip
      \\
      \lambda^{1/2}c_\lambda^{a_m(G)/2} & \ \mbox{if }
      \textcolor{black}{1\gtrsim}c_\lambda
    \gtrsim
    \lambda^{-(r-1)/(e(G)-a_m(G))}, 
    \end{array}
    \right.
  \end{eqnarray}
  \textcolor{black}{as $\lambda$ tends to infinity.}
    In particular,      
    when $G$ has no endpoints $(m=0)$, we have\textcolor{black}{, as $\lambda$
        tends to infinity}, 
  \begin{eqnarray}
\nonumber 
    \Delta_\lambda\asymp\left\{\begin{array}{ll}
    \lambda^{v(G)/2}c_\lambda^{e(G)/2} & \ \mbox{if }\
     c_\lambda
      \lesssim \lambda^{-(v(G)-1)/e(G)}
    ,
     \medskip
      \\
      \displaystyle
      \lambda^{1/2} & \ \mbox{if }\
       \textcolor{black}{1\gtrsim}
      c_\lambda
    \gtrsim
    \lambda^{-(v(G)-1)/e(G)}, 
    \end{array}
    \right. 
  \end{eqnarray}
 and \eqref{mbalanced} becomes the strong balance condition. 
  \end{theorem} 
 Theorem~\ref{fdhskfs}
 extends Corollaries~6.4 and 6.6 of \cite{LiuPrivault}
 without restriction to the dilute and sparse regimes
 considered therein. 
 When $G$ is a tree with $v(G)=r$ vertices and no endpoints
$(m=0)$,
 Theorem~\ref{fdhskfs} yields
  \begin{eqnarray}
 \nonumber
       \Delta_\lambda\asymp\left\{\begin{array}{ll}
          \lambda^{v(G)/2}c_\lambda^{(v(G)-1)/2} & \ \mbox{if }\
      \displaystyle
            c_\lambda
      \lesssim \frac{1}{\lambda},
    \medskip
      \\
      \displaystyle
\lambda^{1/2} & \ \mbox{if }\
    \displaystyle
    \textcolor{black}{1\gtrsim}c_\lambda
    \gtrsim
    \frac{1}{\lambda}      , 
    \end{array}
    \right.
  \end{eqnarray}
  which recovers Corollaries~6.4 and 6.6-$1)$ of \cite{LiuPrivault}
  as particular cases. 
\begin{corollary}
  \label{c37}
   (Normal approximation).
        Let $G$ be a connected graph with $V_G=[r+m]$ for $r\geq 2$ and $m\geq 0$, suppose that
        Assumptions~\ref{assm5-0}-\eqref{assm5-01-2}-\eqref{assm5} are satisfied and the balance condition~\eqref{mbalanced} holds.
        In addition, 
    assume that 
    \begin{equation} 
   \label{fjkl243}
  \textcolor{black}{1\gtrsim}c_\lambda \gg  \lambda^{-\min ( r/e(G) , 1 / a_m(G))}. 
   \end{equation} 
  Then, the normalized subgraph count $\widebar{N}_G$
  converges in distribution to a standard normal random variable \textcolor{black}{as $\lambda$  tends to infinity.}
  In particular, when $m=0$ or $m=1$,
  Condition~\eqref{fjkl243} reduces to
  $$
   \textcolor{black}{1\gtrsim}c_\lambda \gg \lambda^{-r/e(G)},
  $$
  i.e. $\textcolor{black}{1\gtrsim}c_\lambda \gg \lambda^{-v(G)/e(G)}$
  when $G$ has no endpoints $(m=0)$. 
\end{corollary}
\begin{Proof}
  It suffices to note that under \eqref{fjkl243},
  in both cases 
  \begin{enumerate}[i)]
  \item
    $\textcolor{black}{1\gtrsim}c_\lambda \gg  \lambda^{-r/e(G)}$
       if $a_m(G)r < e(G)$, 
   and 
      \item
$    \textcolor{black}{1\gtrsim}c_\lambda
    \gg 
    \lambda^{-1/a_m(G)}$
        if $a_m(G)r \geq e(G)$, 
  \end{enumerate}
  we have $\lim_{\lambda \to \infty} \Delta_\lambda = \infty$
  in Theorem~\ref{fdhskfs},
  and to apply Theorem~1 in \cite{Janson1988}. 
\end{Proof} 
\textcolor{black}{When $a_m(G)=0$, and in particular if $m=0$,
  we have $1/a_m(G)=\infty$, and therefore
  $$\min ( r/e(G) , 1 / a_m(G))=r/e(G)$$
  in \eqref{fjkl243}.} 
 When $c_\lambda \asymp \lambda^{-\alpha}$,
 the normal approximation result
 of Corollary~\ref{c37} holds 
 provided that 
 \begin{equation}
   \label{fjkld34} 
 0 < \alpha <\min \left( \frac{r}{e(G)} , \frac{1}{a_m(G)}\right).
\end{equation} 
 When $m=0$ or $m=1$,
 Condition~\eqref{fjkld34}
 is equivalent to 
 $$
 0 < \alpha < \frac{r}{e(G)}, 
 $$ 
 which also reads
 $$
 0 < \alpha < \frac{r}{e(G)} = \frac{v(G)}{e(G)} 
 $$ 
 in the absence of endpoints $(m=0)$. 

 \medskip

 Theorem~\ref{thmstatu} follows from Theorem~\ref{fdhskfs}
  and the ``main lemmas'' in Chapter~2 of
  \cite{saulis} and in \cite{doering}.
   When $m=0$, the Kolmogorov rate in \eqref{normalapp}
   with $c_\lambda
    \gtrsim
    \lambda^{-(v(G)-1)/e(G)}$ is
 consistent with the rate in 
 Corollary~4.6 of \cite{PS2}
 in the Erd{\H o}s-R\'enyi model, up to an additional
 power $1/(2r-1)$. 
\begin{theorem}
\label{thmstatu}
   (Normal approximation).
  Let $G$ be a connected graph with $V_G=[r+m]$ for $r\geq 2$ and $m\geq 0$.
  Suppose that
  Assumptions~\ref{assm5-0}-\eqref{assm5-01-2}-\eqref{assm5} are satisfied,
  that the balance condition~\eqref{mbalanced} holds,
  that 
   \begin{equation*} 
  \textcolor{black}{1\gtrsim}c_\lambda \gg  \lambda^{-\min ( r/e(G) , 1 / a_m(G))}, 
   \end{equation*} 
   and let $\Delta_\lambda$ be defined in \eqref{d1}. 
\begin{enumerate}[i)]
\item (Kolmogorov bound,
 \cite[Corollary~2.1]{saulis} and \cite[Theorem~2.4]{doering})
  One has
  \begin{equation}
\label{normalapp}
\sup_{x\in\R}|\IP_\lambda (\widebar{N}_G \leq x)-\Phi(x)|\leq \frac{C}{(\Delta_\lambda)^{1/(2r-1)}},
\end{equation}
 where $C>0$ is a constant depending only on $r \geq 2$.
\item (Moderate deviation principle,
  \cite[Theorem~1.1]{doring} and \cite[Theorem~3.1]{doering}).
  Let $( a_\lambda )_{\lambda > 0}$ be a \textcolor{black}{function of $\lambda$} tending to infinity \textcolor{black}{as $\lambda$ tends to infinity}, and such that 
  $$
  \lim_{\lambda \to \infty}
  \frac{a_\lambda}{(\Delta_\lambda)^{1/(2r-1 )}}
  = 0.
  $$
    Then, $ (a_\lambda^{-1}\widebar{N}_G )_{\lambda >0}$ satisfies a moderate deviation principle with speed $a_\lambda^2$ and rate function $x^2/2$. 
\item (Concentration inequality,
  corollary of \cite[Lemma~2.4]{saulis} and \cite[Theorem~2.5]{doering}).
  For any $x\geq 0$ and sufficiently large $\lambda$, 
  \begin{equation}
    \nonumber 
\IP_\lambda (|\widebar{N}_G |\geq x)\le2\exp\left(-\frac14\min\left( \frac{x^2}{2^{r}},(x\Delta_\lambda)^{1/r}\right) \right).
\end{equation} 
\item (Normal approximation with Cram\'er corrections, \cite[Lemma~2.3]{saulis} \textcolor{black}{and \cite[Theorem~2.3]{doering}}). There exists a constant $c>0$ such that for all $\lambda \geq 1$ and $x\in(0,c(\Delta_\lambda)^{1/(2r-1 )})$ we have 
  $$
  \frac{\IP_\lambda (\widebar{N}_G \geq x)}{1-\Phi(x)}= \left(1+O\left(\frac{x+1}{(\Delta_\lambda)^{1/(2r-1 )}}\right)\right) \exp \big( \widetilde{L}(x) \big)
  $$
  and
  $$
  \frac{\IP_\lambda (\widebar{N}_G \leq -x)}{\Phi(-x)}=\left(1+O\left(\frac{x+1}{(\Delta_\lambda)^{1/(2r -1 )}}\right)\right) \exp \big( \widetilde{L}(-x) \big), 
$$ 
  where
  \textcolor{black}{ 
    $\widetilde{L}(x):= (x/c)^3(\Delta_\lambda)^{-3/(2r-1 )}\theta$, 
    for some $\theta\in[-1,1]$ depending on
$x\in(0,c(\Delta_\lambda)^{1/(2r-1 )})$.}
    \end{enumerate}
\end{theorem}
  \subsubsection*{Poisson approximation}
\noindent
In what follows, 
$|\mathrm{Aut}_{\bullet}(G)|$ stands for the number of automorphisms of \textcolor{black}{the} induced subgraph $H\subset G$ with $V_H=[r]$. \textcolor{black}{For example, taking $G$ in Figure~\ref{fig:graph1}, 
  the induced graph of $G$ is a path $H$ with vertex set $[4]$ and
  edge set $E_H=\{\{1,2\},\{1,3\},\{3,4\}\}$, which gives $|\mathrm{Aut}_{\bullet}(G)|=2$.}
  \textcolor{black}{The} condition $a_m(G) r \textcolor{black}{\le} e(G)$ in the
  Poisson limit Theorem~\ref{poslim}
  always holds when $m=0$, 
  and it holds for \textcolor{black}{trees}, cycles and complete graphs 
  when $m=1$. 
  \begin{theorem}\label{poslim}
    (Poisson approximation).
    Let $G$ be a connected graph with $V_G=[r+m]$
    for $r\geq 2$ and $m\geq 0$,
    suppose that 
    Assumptions~\ref{assm5-0}-\eqref{assm5-01-2}-\eqref{assm5-1-2}
    are satisfied. 
   If the balance condition~\eqref{mbalanced} holds together with
  $a_m(G) r \leq e(G)$ and 
$$
    \lim_{\lambda \to \infty} \lambda c_\lambda^{e(G)/r} =c>0, 
    $$ 
  then the subgraph count $\widehat{N}_G:=N_G/|\mathrm{Aut}_{\bullet}(G)|$
  converges in distribution to a Poisson random variable with mean
 \begin{equation}
\nonumber 
  \mu_\varphi:=
  \frac{c^r}{|\mathrm{Aut}_{\bullet}(G)|}\int_{(\R^d)^r}\prod_{\textcolor{black}{\{i,j\}}\in E_G}\varphi(x_i,x_j)
  \ \! \mu(\mathrm{d}x_1)\cdots\mu(\mathrm{d}x_r),
  \end{equation}
 where $x_{r+i}:=y_i$ for $i=1,\dots,m$.
\end{theorem}
 \textcolor{black}{The proof of Theorem~\ref{poslim} is postponed to the end
  of Section~\ref{cumulantrates}.}

    \medskip 
  
 \noindent
 By the first and second moment methods,
 see \cite[Page~54]{JLR}, we have 
\begin{equation}
\label{firstsecm}
\frac{(\E_\lambda [X])^2}{\E_\lambda [X^2]}\le\IP_\lambda (X>0)\le \E_\lambda [X]
\end{equation}
 for any non-negative integer-valued random variable $X$,
 which yields \textcolor{black}{the}
  following \textcolor{black}{threshold} result for subgraph containment
 as a consequence of Corollary~\ref{themo-1-1}. 
\begin{theorem}\label{pt:contain}
  (\textcolor{black}{Threshold} for subgraph containment).
  Let $G$ be a connected graph with $V_G=[r+m]$ for $r\geq 2$ and $m\geq 0$,
  and suppose that
    Assumptions~\ref{assm5-0}-\eqref{assm5-01-2}-\eqref{assm5-1-2}
    are satisfied. 
  If the balance condition~\eqref{mbalanced} holds together with
  $a_m(G) r \leq e(G)$, 
  then we have the \textcolor{black}{following threshold} results: 
\begin{enumerate}[a)] 
\item
  \label{part-a}
  $\lim_{\lambda \to \infty} \IP_\lambda (N_G = 0) = 1$ if $c_\lambda\ll\lambda^{-r/e(G)}$,
\item $\lim_{\lambda \to \infty} \IP_\lambda (N_G = 0) =
  \textcolor{black}{e^{- \nu_\varphi}}$ if $c_\lambda\sim \lambda^{-r/e(G)}$,
  with
  $$ 
  \nu_\varphi:=
  \frac{1}{|\mathrm{Aut}_{\bullet}(G)|}\int_{(\R^d)^r}\prod_{\textcolor{black}{\{i,j\}}\in E_G}\varphi(x_i,x_j)
  \ \! \mu(\mathrm{d}x_1)\cdots\mu(\mathrm{d}x_r), 
$$   
  where $x_{r+i}:=y_i$ for $i=1,\dots,m$, 
\item $\lim_{\lambda \to \infty} \IP_\lambda (N_G = 0) = 0$ if
 $1 \gtrsim c_\lambda\gg\lambda^{-r/e(G)}$. 
\end{enumerate} 
\end{theorem}
\begin{Proof}
$(a)$ Since $\E_\lambda [N_G]\asymp\lambda^r c_\lambda^{e(G)}$, if $c_\lambda\ll\lambda^{-r/e(G)}$,
 we know that $\lim_{\lambda \to \infty} \E_\lambda [N_G] = 0$,
 and we conclude by the first moment method in \eqref{firstsecm}. 
 \medskip

\noindent 
 $(b)$ 
Since $c_\lambda\sim \lambda^{-r/e(G)}$,
we have 
 $c = \lim_{\lambda \to \infty} \lambda c_\lambda^{e(G)/r} =1$,
and we conclude by Theorem~\ref{poslim}.
 
 \medskip

\noindent 
 $(c)$ 
 From Corollary~\ref{themo-1-1} we know that if
  $c_\lambda\gg\lambda^{-(r-1)/(e(G)-a_m(G))}$, \textcolor{black}{then}
$$\kappa_2(N_G)\asymp\lambda^{2r-1}c_{\lambda}^{2e(G)-a_m(G)},$$
and if $c_\lambda\ll\lambda^{- ( r-1 ) / ( e(G)-a_m(G) ) }$, \textcolor{black}{then}
$$\kappa_2(N_G)\asymp\lambda^r c_{\lambda}^{e(G)}.
$$
If $a_m(G) \leq e(G)/r$ and $c_\lambda\gg\lambda^{-r/e(G)}$, \textcolor{black}{then} we have 
\begin{equation}
\frac{(\E_\lambda [N_G])^2}{\E_\lambda [(N_G)^2]}\asymp\frac{\lambda^{2r}c_\lambda^{2e(G)}}{\lambda^{2r}c_\lambda^{2e(G)}+\kappa_2(N_G)}\asymp\frac1{1+\big(\lambda^r c_{\lambda}^{e(G)}\big)^{-1}\vee \big(\lambda c_\lambda^{a_m(G)}\big)^{-1}}\to 1 ,
\end{equation}
 and we conclude by the second moment method in \eqref{firstsecm}. 
\end{Proof}
\section{Diagram representation of cumulants} 
\label{sec:diag}
\noindent 
This section introduces the diagram framework used for the
expansion of cumulants as sums over partitions
in Proposition~\ref{djkdsf34} below. 

\medskip

We start with basic notation on set partitions, see e.g. \cite{peccatitaqqu}, \cite{LiuPrivault}.
For any finite set $b$, we let $\Pi(b)$ denote the collection of set partitions of $b$. For two set partitions $\rho_1,\rho_2\in\Pi(b)$, we say $\rho_1$ is
coarser than $\rho_2$
 (i.e. $\rho_2$ is finer than $\rho_1$), 
 and we write it as $\rho_2\preceq\rho_1$, 
if and only if each block of $\rho_2$ is contained in a block of $\rho_1$.
We use $\rho_1\vee\rho_2$ for the
finest partition which is coarser than both $\rho_1$ and $\rho_2$, and
denote by $\rho_1\wedge\rho_2$ the coarsest partition
which is finer than both of $\rho_1$ and $\rho_2$. We also let $\widehat{1}:=\{b\}$ 
denote the coarsest partition of $b$,
 whereas $\widehat{0}$
 stands for the partition made of singletons.
 \begin{definition}
   Given $r\geq 2$ and $n\geq 1$ 
 we let $\pi$ denote the partition
 $\pi:=\{\pi_1,\dots,\pi_n\}\in\Pi([n]\times[r])$
 of $[n]\times[r]$ defined as 
 $$\pi_i:=\{(i,j) \ : \ 1\le j\le r\},
 \quad i =1,\ldots , n.
 $$
 \begin{enumerate}[a)]
 \item
   A partition $\rho\in\Pi([n]\times[r])$ is said to be {\it non-flat} if $\rho\wedge\pi=\widehat{0}$, and {\it connected} if $\rho\vee\pi=\widehat{1}$.
 \item
   We let $\Pi_{\widehat{1}}([n]\times[r])$ denote the collection of all connected partitions of $[n]\times[r]$, and let 
      $$
      {\rm CNF}(n,r) :=\big\{\rho \ : \ \rho\in\Pi_{\widehat{1}}( [n]\times [r]),\ \rho\wedge\pi=\widehat{0}\big\}
      $$
 denote     the set of all connected and non-flat partitions of $[n]\times [r]$, for $n, r\geq 1$.
 \end{enumerate}
 \end{definition}
\begin{example}
 Figure~\ref{fig:diagram0} presents an example
 of non-flat connected partition in ${\rm CNF}(3,4)$. 

\begin{figure}[H]
\centering
 \begin{tikzpicture}[scale=0.8] 
 \draw[black, thick] (0,0) rectangle (5,4);
 \node[anchor=east,font=\small] at (0.8,3) {1};
 \node[anchor=east,font=\small] at (0.8,2) {2};
 \node[anchor=east,font=\small] at (0.8,1) {3};
 \node[anchor=south,font=\small] at (1,0) {1};
 \node[anchor=south,font=\small] at (2,0) {2};
 \node[anchor=south,font=\small] at (3,0) {3};
 \node[anchor=south,font=\small] at (4,0) {4};
 
 \filldraw [black] (1,1) circle (2pt);
 \filldraw [black] (2,1) circle (2pt);
 \filldraw [black] (3,1) circle (2pt);
 \filldraw [black] (4,1) circle (2pt);
 \filldraw [black] (1,2) circle (2pt);
 \filldraw [black] (2,2) circle (2pt);
 \filldraw [black] (3,2) circle (2pt);
 \filldraw [black] (4,2) circle (2pt);
 \filldraw [black] (1,3) circle (2pt);
 \filldraw [black] (2,3) circle (2pt);
 \filldraw [black] (3,3) circle (2pt);
 \filldraw [black] (4,3) circle (2pt);
 
 \draw[very thick] (1,3) -- (1,2) -- (1,1);
  \draw[very thick] (3,3) -- (3,2);
 \draw[very thick] (2,2) -- (2,1);
 \draw[very thick] (3,1) -- (4,2);
 
 \begin{pgfonlayer}{background}
    \filldraw [line width=4mm,black!3]
      (0.2,0.2)  rectangle (4.8,3.8);
  \end{pgfonlayer}
 \end{tikzpicture}
 \caption{Non-flat connected partition of
  $[3]\times[4]$. }
 \label{fig:diagram0}
 \end{figure}
\end{example}
\noindent
 We also note the following lemma. 
\begin{lemma}\cite[Lemma~2.5]{LiuPrivault} 
  \label{restrict-partition}
  Let $n\geq 2$.
  For any connected partition 
  $\rho\in \Pi_{\widehat{1}}( [n]\times [r])$
  there exists $i\in \{1,\dots,n \}$
  such that the set partition 
    $\{b\backslash\pi_i:b\in\rho\}$ 
  of $\{1,\dots,i-1,i+1,\dots,n \}\times [r]$
  is connected.
\end{lemma}
\noindent
In \cite{LiuPrivault}, a graphical diagram language has been designed for the cumulant representation of subgraph \textcolor{black}{counts} in the random-connection model, and extended to the case of subgraphs
containing fixed endpoints in c.f. \cite{LiuPri23b}.
 \begin{definition}\label{defgraph-1}
  {\cite{LiuPri23b}}
  \textcolor{black}{Let $G$ be a connected graph of order $r+m$ satisfying Assumption~\ref{a0}, and}
  let $\rho=\{b_1,\dots,b_{|\rho |} \} \in \Pi([n]\times[r])$, $n \geq 1$, be
  a partition of $[n]\times[r]$.
    We let $\widebar{\rho}_G$ denote the connected multigraph
  built on \textcolor{black}{ $[m]\cup([n]\times[r])$}, which is constructed as follows.  
     \begin{enumerate}  
     \item \textcolor{black}{For} all $i\in[n]$ and $j_1,j_2\in[r]$, $j_1\ne j_2$, an edge links $(i,j_1)$ to $(i,j_2)$ iff $(j_1,j_2)\in E_G$;
     \item \textcolor{black}{For} all $k\in[m]$, $i\in[n]$ and $j\in[r]$, an edge links $(k)$ to $(i,j)$ iff $(r+k,j)\in E_G$;
     \item \textcolor{black}{For} all $i\in[|\rho |]$, \textcolor{black}{all elements in the same block $b_i$ are regarded as one vertex.}
    \end{enumerate}
 In addition, we let $\rho_G$ be the graph constructed
from the multigraph $\widebar{\rho}_G$ 
by \textcolor{black}{replacing multiple edges with simple edges}
in $\widebar{\rho}_G$. 
\end{definition}
In what follows, the blocks of any
given partition $\rho=\{b_1,\dots,b_{|\rho |} \}$
in $\Pi([n]\times[r])$ are ordered along the lexicographic order
on $[n]\times[r]$,
 by ordering the blocks according to their smallest elements.
From the above construction, the vertices of $\rho_G$ originate
from terminal nodes $[m]$ and blocks $b_1,\ldots , b_{|\rho |}$ of $\rho$.
We can further denote the vertex set of the graph $\rho_G$ as $V(\rho_G):=\left[|\rho|+m\right]$ according to the rule that \textcolor{black}{the $m$ terminal nodes} follows $b_1,\ldots , b_{|\rho |}$ in order.

\textcolor{black}{See also \cite{khorunzhiy} for a diagram representation used
 for lines and cycles in the Erd{\H o}s-R\'enyi model, and \cite{FGY23} for a graphical representation defined for the $U$-statistics of determinantal point processes.}
\begin{example}
  Consider $\rho\in\Pi([3]\times[4])$ as
  in Figure~\ref{fig:diagram0}, with 
  \begin{eqnarray*}
    \rho&=&\{\{(1,1),(2,1),(3,1)\},
    \{(1,2)\},\{(1,3),(2,3)\},\\
    &&\ \ \ \{(1,4)\},\{(2,2),(3,2)\},\{(2,4),(3,3)\},\{(3,4)\}\}, 
  \end{eqnarray*}
  and let $G$ be the connected graph with vertex set $V_G=[5]$,
  represented in Figure~\ref{fig:graph1}.
\begin{figure}[H]
  \centering
   \begin{tikzpicture}[scale=0.8] 
\draw[black, thick] (0,0) rectangle (6,2);
\filldraw [red] (1,1) circle (2pt);
\foreach \i in {2,3,4,5}
{
\filldraw [blue] (\i,1) circle (2pt);
}
\draw[thick] (1,1) -- (3,1);
\draw[thick] (4,1) -- (5,1);
\draw[thick] (2,1) .. controls (3,0.5) .. (4,1);
\node[anchor=south,font=\small] at (2,1) {1};
\node[anchor=south,font=\small] at (3,1) {2};
\node[anchor=south,font=\small] at (4,1) {3};
\node[anchor=west,font=\small] at (5,1) {4};
\node[anchor=east,red,font=\small] at (1,1) {5};
  \begin{pgfonlayer}{background}
    \filldraw [line width=4mm,black!3]
      (0.2,0.2)  rectangle (5.8,1.8);
  \end{pgfonlayer}
  \end{tikzpicture}
   \caption{Connected graph $G$ on five vertices including one endpoint, with $r=4$ and $m=1$.}
   \label{fig:graph1}
   \end{figure}
   
 \vspace{-0.4cm}

\noindent
  Figure~\ref{fig:diagram1} presents the multigraph $\widebar{\rho}_G$ and
  corresponding graph $\rho_G$. 

  \smallskip

\begin{figure}[H]
\captionsetup[subfigure]{font=footnotesize}
\centering
\subcaptionbox{Multigraph $\widebar{\rho}_G$ before merging edges and vertices.}[.49\textwidth]{
\begin{tikzpicture}
\draw[step=1cm, very thin, gray!40] (0,0) grid (6,4);
\draw[black, thick] (0,0) rectangle (6,4);
\filldraw [red] (1,2) circle (2pt);
\foreach \i in {2,3}
{
\filldraw [black] (2,\i) circle (2pt);
\filldraw [black] (3,\i) circle (2pt);
\filldraw [black] (4,\i) circle (2pt);
\filldraw [black] (5,\i) circle (2pt);
\draw[thick, dash dot,blue] (2,\i) -- (3,\i);
\draw[thick, dash dot,blue] (2,\i) .. controls (3,\i-0.5) .. (4,\i);
\draw[thick, dash dot,blue] (1,2) -- (2,\i);
\draw[thick, dash dot,blue] (4,\i) -- (5,\i);
}

\filldraw [black] (2,1) circle (2pt);
\filldraw [black] (3,1) circle (2pt);
\filldraw [black] (4,1) circle (2pt);
\filldraw [black] (5,1) circle (2pt);
\draw[thick, dash dot,blue] (2,1) -- (3,1);
\draw[thick, dash dot,blue] (2,1) .. controls (3,1+0.5) .. (4,1);
\draw[thick, dash dot,blue] (1,2) -- (2,1);
\draw[thick, dash dot,blue] (4,1) -- (5,1);

\node[anchor=north,font=\scriptsize] at (2,1) {(3,1)};
\node[anchor=north,font=\scriptsize] at (3,1) {(3,2)};
\node[anchor=north,font=\scriptsize] at (4,1) {(3,3)};
\node[anchor=north,font=\scriptsize] at (5,1) {(3,4)};
\node[anchor=south,font=\scriptsize] at (2,3) {(1,1)};
\node[anchor=south,font=\scriptsize] at (3,3) {(1,2)};
\node[anchor=south,font=\scriptsize] at (4,3) {(1,3)};
\node[anchor=south,font=\scriptsize] at (5,3) {(1,4)};
\node[anchor=south,font=\scriptsize] at (2,2) {(2,1)};
\node[anchor=south,font=\scriptsize] at (3,2) {(2,2)};
\node[anchor=south,font=\scriptsize] at (4,2) {(2,3)};
\node[anchor=south,font=\scriptsize] at (5,2) {(2,4)};

\node[anchor=east,font=\small] at (1,2) {$1$};
\draw[very thick] (3,1) -- (3,2);
\draw[very thick] (4,1) -- (5,2);
\draw[very thick] (2,3) -- (2,1);
\draw[very thick] (4,2) -- (4,3);
  \begin{pgfonlayer}{background}
    \filldraw [line width=4mm,black!3]
      (0.2,0.2)  rectangle (5.8,3.8);
  \end{pgfonlayer}
\end{tikzpicture}}
\subcaptionbox{Graph $\rho_G$ after merging edges and vertices.}[.49\textwidth]{
\begin{tikzpicture}
\draw[step=1cm, very thin, gray!40] (0,0) grid (6,4);
  \draw[black, thick] (0,0) rectangle (6,4);
  \filldraw [red] (1,2) circle (2pt);
  \filldraw [black] (2,2) circle (2pt);
  \filldraw [black] (3,1) circle (2pt);
  \filldraw [black] (3,3) circle (2pt);
  \filldraw [black] (5,1) circle (2pt);
  \filldraw [black] (4,2) circle (2pt);
  \filldraw [black] (5,3) circle (2pt);
  \filldraw [black] (4,1) circle (2pt);
  \draw[thick,blue] (1,2) -- (2,2) -- (4,2);
  \draw[thick,blue] (2,2) -- (3,3);
  \draw[thick,blue] (2,2) -- (3,1);
  \draw[thick,blue] (5,3) -- (4,2) -- (4,1);
  
  \draw[thick,blue] (2,2) -- (4,1) -- (5,1);

  \node[anchor=east,font=\footnotesize] at (1,2) {8};
  \node[anchor=south,font=\footnotesize] at (3,3) {2};
  \node[anchor=north,font=\footnotesize] at (3,1) {5};
  \node[anchor=south,font=\footnotesize] at (5,3) {4};
  \node[anchor=south,font=\footnotesize] at (2,2) {1};
  \node[anchor=south,font=\footnotesize] at (4,2) {3};
  
  \node[anchor=north,font=\footnotesize] at (5,1) {7};
  \node[anchor=north,font=\footnotesize] at (4,1) {6};

  \begin{pgfonlayer}{background}
    \filldraw [line width=4mm,black!3]
      (0.2,0.2)  rectangle (5.8,3.8);
  \end{pgfonlayer}
\end{tikzpicture}}
\caption{
  Example of graph $\rho_G$
   with $n=3$, $r=4$, and $m=1$.}
\label{fig:diagram1}
\end{figure}

\vspace{-0.5cm}

\end{example}
The above framework allows
us to state moment and cumulant formulas
 for the subgraph counts $N_G$.
\begin{definition}
  For $\rho=\{b_1,\dots,b_{|\rho|}\}\in\Pi([n]\times[r])$ and $j\in[m]$, we denote by 
  \begin{equation}
    {\cal A}^\rho_j:=\{k\in[|\rho|]:~\exists (s,i)\in b_k~\mathrm{s.t.}~\{i,r+j\}\in E_G\}
  \end{equation}
the neighborhood of the vertex $(|\rho|+j)$ in the graph $\rho_G$.
\end{definition}
\noindent
From \cite{LiuPri23b}, we have the following moment and cumulant representation
for $N_G$.
\begin{prop}
  \label{djkdsf34}
  Let $n\geq 1$\textcolor{black}{. Then}, the $n$-th moments and $n$-th cumulants of $N_G$ admit the expressions
  \begin{equation}
  \label{cum-eq} 
    \E_\lambda [(N_G)^n]=\sum_{\rho \in \Pi( [n] \times [r])
    \atop
    {\rho \wedge \pi = \widehat{0}
    \atop {\rm (non-flat)}
    }
    }
    F_\lambda (\rho)
    \quad
    \mbox{and}
    \quad 
  \kappa_n(N_G)=
  \sum_{\rho \in \Pi_{\widehat{1}} ( [n] \times [r])
    \atop
    {\rho \wedge \pi = \widehat{0}
    \atop {\rm (non-flat \ \! connected)}
    }
    }
F_\lambda (\rho),
\end{equation} 
 where $F_\lambda (\rho)$, $\rho \in \Pi( [n] \times [r])$, is defined as 
\begin{equation}
\label{fjhkldf231} 
F_\lambda (\rho):=
\lambda^{|\rho |}
\int_{(\R^d)^{|\rho|}}\prod_{\substack{ 
  1 \leq j \leq m
  \\ i\in {\cal A}^\rho_j}}
\varphi_\lambda(x_i,y_j)
\ \prod_{
  \substack{1 \leq k \textcolor{black}{<} l\le|\rho|
    \\
    \{ k , l \}\in E(\rho_G) 
}}\varphi_\lambda(x_k,x_l) \mu ( \mathrm{d}x_1 ) \cdots \mu ( \mathrm{d}x_{|\rho|}). 
\end{equation}
\end{prop}
\section{Planar representation of partition diagrams} 
\label{plot1}
\noindent
In this section, we introduce a planar representation
that will allow us to determine the leading 
 partition diagrams in the
moment and cumulant expressions of Proposition~\ref{djkdsf34}.
\begin{definition}\label{xyplot}
 Let $G$ be a connected graph with $V_G=[r+m]$, for some $r\geq 2$ and $m\geq 0$.
 For $n\geq 2$, we let 
\begin{equation}
\nonumber
\Sigma_n(G,m):=\big\{
( x ( \rho_G ) , y ( \rho_G ) ) :=
(nr+m-v(\rho_G),ne(G)-e(\rho_G)): \ \rho \in
        {\rm CNF} (n,r)
        \big\}, 
\end{equation}
where, for every partition $\rho\in\Pi([n]\times[r])$,
$\rho_G$ is the graph associated to $\rho$ by Definition~\ref{defgraph-1}. 
 \end{definition}
\noindent
\begin{example}
  Let $G = C_3$ be a triangle with no endpoint,
  i.e. $r=3$ and $m=0$. We have
 $$
 \left\{
 \begin{array}{l}
  \Sigma_2(C_3,0) = \{(3,3),(2,1),(1,0)\},
   \medskip
   \\ 
  \Sigma_3(C_3,0) = \{(6, 6), (5, 4), (4, 3), (5, 3), (4, 2), (4, 1), (3, 1), (3, 0), (2, 0)\},
     \medskip
   \\ 
   \Sigma_4(C_3,0) = \{(9, 9), (8, 7), (7, 6), (8, 6), (7, 5), (7, 4), (6, 4), (6, 3),(5, 3), (7, 3),
   \smallskip
   \\
   \qquad \qquad \quad \ \ \ \! \!
    (6, 2), (5, 2), (7, 2), (6, 1), (5, 1), (4, 1), (6, 0), (5, 0), (4, 0), (3, 0)\},
 \end{array}
 \right.
   $$ 
 see Figure~\ref{fig1}.
 
  \vspace{-0.4cm}
 
\begin{figure}[H]
  \begin{subfigure}{.5\textwidth}
    \vskip1cm
  \centering
\begin{tikzpicture}[scale=0.8] 
\draw[step={(1cm,0.5cm)}, very thin, gray!40] (0,0) grid (7,3.5);
\draw[->] (0,0) -- (0,3.5);
\node[anchor=east,font=\small] at (-0.05,3.5) {$y$};
\draw[->] (0,0) -- (7,0);
\node[anchor=north,font=\small] at (7,-0.05) {$x$};
\draw[thick] (-0.05,0.5) -- (0.05,0.5);
\node[anchor=east,font=\small] at (-0.05,0.5) {1};
\draw[thick] (-0.05,1) -- (0.05,1);
\node[anchor=east,font=\small] at (-0.05,1) {2};
\draw[thick] (-0.05,2) -- (0.05,2);
\node[anchor=east,font=\small] at (-0.05,2) {4};
\draw[thick] (-0.05,3) -- (0.05,3);
\node[anchor=east,font=\small] at (-0.05,3) {6};
\draw[thick] (-0.05,1.5) -- (0.05,1.5);
\node[anchor=east,font=\small] at (-0.05,1.5) {3};
\draw[thick] (-0.05,2.5) -- (0.05,2.5);
\node[anchor=east,font=\small] at (-0.05,2.5) {5};
\node[anchor=east,font=\small] at (-0.05,3) {6};
\draw[thick] (1,-0.05) -- (1,0.05);
\node[anchor=north,font=\small] at (1,-0.05) {1};
\draw[thick] (2,-0.05) -- (2,0.05);
\node[anchor=north,font=\small] at (2,-0.05) {2};
\draw[thick] (3,-0.05) -- (3,0.05);
\node[anchor=north,font=\small] at (3,-0.05) {3};
\draw[thick] (4,-0.05) -- (4,0.05);
\node[anchor=north,font=\small] at (4,-0.05) {4};
\draw[thick] (5,-0.05) -- (5,0.05);
\node[anchor=north,font=\small] at (5,-0.05) {5};
\draw[thick] (6,-0.05) -- (6,0.05);
\node[anchor=north,font=\small] at (6,-0.05) {6};
\draw[very thick,red] (2,0) -- (6,3);
\filldraw [blue] (2,0) circle (2.0pt);
\filldraw [blue] (3,0) circle (2.0pt);
\filldraw [blue] (3,0.5) circle (2.0pt);
\filldraw [blue] (4,0.5) circle (2.0pt);
\filldraw [blue] (4,1) circle (2.0pt);
\filldraw [blue] (5,1.5) circle (2.0pt);
\filldraw [blue] (4,1.5) circle (2.0pt);
\filldraw [blue] (5,2) circle (2.0pt);
\filldraw [blue] (6,3) circle (2.0pt);
    \end{tikzpicture}
\caption{Subgraph plot of $\Sigma_3(C_3,0)$.} 
\end{subfigure}
\begin{subfigure}{.5\textwidth}
  \centering
  \begin{tikzpicture}[scale=0.8] 
\draw[step={(0.5cm,0.5cm)}, very thin, gray!40] (0,0) grid (5,5);
\draw[->] (0,0) -- (0,5);
\node[anchor=east,font=\small] at (-0.05,5) {$y$};
\draw[->] (0,0) -- (5,0);
\node[anchor=north,font=\small] at (5,-0.05) {$x$};
\draw[thick] (-0.05,0.5) -- (0.05,0.5);
\node[anchor=east,font=\small] at (-0.05,0.5) {1};
\draw[thick] (-0.05,1) -- (0.05,1);
\node[anchor=east,font=\small] at (-0.05,1) {2};
\draw[thick] (-0.05,1.5) -- (0.05,1.5);
\node[anchor=east,font=\small] at (-0.05,1.5) {3};
\draw[thick] (-0.05,2) -- (0.05,2);
\node[anchor=east,font=\small] at (-0.05,2) {4};
\draw[thick] (-0.05,2.5) -- (0.05,2.5);
\node[anchor=east,font=\small] at (-0.05,2.5) {5};
\draw[thick] (-0.05,3) -- (0.05,3);
\node[anchor=east,font=\small] at (-0.05,3) {6};
\draw[thick] (-0.05,3.5) -- (0.05,3.5);
\node[anchor=east,font=\small] at (-0.05,3.5) {7};
\draw[thick] (-0.05,4) -- (0.05,4);
\node[anchor=east,font=\small] at (-0.05,4) {8};
\draw[thick] (-0.05,4.5) -- (0.05,4.5);
\node[anchor=east,font=\small] at (-0.05,4.5) {9};
\draw[thick] (0.5,-0.05) -- (0.5,0.05);
\node[anchor=north,font=\small] at (0.5,-0.05) {1};
\draw[thick] (1,-0.05) -- (1,0.05);
\node[anchor=north,font=\small] at (1,-0.05) {2};
\draw[thick] (1.5,-0.05) -- (1.5,0.05);
\node[anchor=north,font=\small] at (1.5,-0.05) {3};
\draw[thick] (2,-0.05) -- (2,0.05);
\node[anchor=north,font=\small] at (2,-0.05) {4};
\draw[thick] (2.5,-0.05) -- (2.5,0.05);
\node[anchor=north,font=\small] at (2.5,-0.05) {5};
\draw[thick] (3,-0.05) -- (3,0.05);
\node[anchor=north,font=\small] at (3,-0.05) {6};
\draw[thick] (3.5,-0.05) -- (3.5,0.05);
\node[anchor=north,font=\small] at (3.5,-0.05) {7};
\draw[thick] (4,-0.05) -- (4,0.05);
\node[anchor=north,font=\small] at (4,-0.05) {8};
\draw[thick] (4.5,-0.05) -- (4.5,0.05);
\node[anchor=north,font=\small] at (4.5,-0.05) {9};
\draw[very thick,red] (1.5,0) -- (4.5,4.5);
\filldraw [blue] (1.5,0) circle (2.0pt);
\filldraw [blue] (2,0) circle (2.0pt);
\filldraw [blue] (2.5,0) circle (2.0pt);
\filldraw [blue] (3,0) circle (2.0pt);
\filldraw [blue] (2,0.5) circle (2.0pt);
\filldraw [blue] (2.5,0.5) circle (2.0pt);
\filldraw [blue] (3,0.5) circle (2.0pt);
\filldraw [blue] (3.5,1) circle (2.0pt);
\filldraw [blue] (3.0,1) circle (2.0pt);
\filldraw [blue] (2.5,1) circle (2.0pt);
\filldraw [blue] (3.5,1.5) circle (2.0pt);
\filldraw [blue] (2.5,1.5) circle (2.0pt);
\filldraw [blue] (3,1.5) circle (2.0pt);
\filldraw [blue] (3,2) circle (2.0pt);
\filldraw [blue] (3.5,2) circle (2.0pt);
\filldraw [blue] (3.5,2.5) circle (2.0pt);
\filldraw [blue] (4,3) circle (2.0pt);
\filldraw [blue] (3.5,3) circle (2.0pt);
\filldraw [blue] (4,3.5) circle (2.0pt);
\filldraw [blue] (4.5,4.5) circle (2.0pt);
\end{tikzpicture}
  \caption{Subgraph plot of $\Sigma_4(C_3,0)$.} 
\end{subfigure}
\caption{Set $\Sigma_n(C_3,0)$ and upper boundary of its convex hull (in red) for $n=3,4$.}
\label{fig1}
\end{figure}

\vspace{-0.4cm}

\noindent 
 Figure~\ref{fig1} and the following ones 
 can be plotted after loading the SageMath code presented in \textcolor{black}{the}
 appendix and running the following commands.

\begin{table}[H] 
\centering
\scriptsize 
 \resizebox{0.9\textwidth}{!}
    {
  \begin{tabular}{|ll|l|} 
 \hline
 \multicolumn{2}{|l|}{
   \EscVerb{G = [[1,2],[2,3],[3,1]]; EP = []; SG3=convexhull(3,G,EP); SG4=convexhull(4,G,EP)
   }
 }
 \\
 \hline
 \multicolumn{2}{|l|}{
\EscVerb{Polyhedron(SG3).plot(color = "pink")+point(SG3,color = "blue",size=20)}
}
 \\
 \hline
 \multicolumn{2}{|l|}{
\EscVerb{Polyhedron(SG4).plot(color = "pink")+point(SG4,color = "blue",size=20)}
 }
 \\
\hline
\end{tabular}
}
\end{table} 

\vspace{-0.4cm}
\end{example}
\begin{example}
  Let $G=C_4$ be a $4$-cycle with no endpoint,
  i.e. $r=4$ and $m=0$.
  Here, $\Sigma_2(C_4,0)$ and $\Sigma_3(C_4,0)$
  are plotted in Figure~\ref{fig2}. 

  \vspace{-0.2cm}
 
\begin{figure}[H]
\begin{subfigure}{.5\textwidth}
  \vskip2cm
    \centering
\begin{tikzpicture}[scale=0.8] 
\draw[step={(1cm,0.5cm)}, very thin, gray!40] (0,0) grid (5,2.5);
\draw[->] (0,0) -- (0,2.5);
\node[anchor=east,font=\small] at (-0.05,2.5) {$y$};
\draw[->] (0,0) -- (5,0);
\node[anchor=north,font=\small] at (5,-0.05) {$x$};
\draw[thick] (-0.05,0.5) -- (0.05,0.5);
\node[anchor=east,font=\small] at (-0.05,0.5) {1};
\draw[thick] (-0.05,1) -- (0.05,1);
\node[anchor=east,font=\small] at (-0.05,1) {2};
\draw[thick] (-0.05,1.5) -- (0.05,1.5);
\node[anchor=east,font=\small] at (-0.05,1.5) {3};
\draw[thick] (-0.05,2) -- (0.05,2);
\node[anchor=east,font=\small] at (-0.05,2) {4};
\draw[thick] (1,-0.05) -- (1,0.05);
\node[anchor=north,font=\small] at (1,-0.05) {1};
\draw[thick] (2,-0.05) -- (2,0.05);
\node[anchor=north,font=\small] at (2,-0.05) {2};
\draw[thick] (3,-0.05) -- (3,0.05);
\node[anchor=north,font=\small] at (3,-0.05) {3};
\draw[thick] (4,-0.05) -- (4,0.05);
\node[anchor=north,font=\small] at (4,-0.05) {4};
\draw[very thick,red] (1,0) -- (4,2);
\filldraw [blue] (1,0) circle (2.0pt);
\filldraw [blue] (2,0) circle (2.0pt);
\filldraw [blue] (2,0.5) circle (2.0pt);
\filldraw [blue] (3,0.5) circle (2.0pt);
\filldraw [blue] (3,1) circle (2.0pt);
\filldraw [blue] (4,1) circle (2.0pt);
\filldraw [blue] (4,2) circle (2.0pt);
\end{tikzpicture}
\caption{Subgraph plot of $\Sigma_2(C_4,0)$.} 
\end{subfigure} 
\begin{subfigure}{.5\textwidth}
  \centering
\begin{tikzpicture}[scale=0.8] 
\draw[->] (0,0) -- (0,4.5);
\draw[step={(0.5cm,0.5cm)}, very thin, gray!40] (0,0) grid (4.5,4.5);
\node[anchor=east,font=\small] at (-0.05,4.5) {$y$};
\draw[->] (0,0) -- (5,0);
\node[anchor=north,font=\small] at (4.5,-0.05) {$x$};
\draw[thick] (-0.05,0.5) -- (0.05,0.5);
\node[anchor=east,font=\small] at (-0.05,0.5) {1};
\draw[thick] (-0.05,1) -- (0.05,1);
\node[anchor=east,font=\small] at (-0.05,1) {2};
\draw[thick] (-0.05,1.5) -- (0.05,1.5);
\node[anchor=east,font=\small] at (-0.05,1.5) {3};
\draw[thick] (-0.05,2) -- (0.05,2);
\node[anchor=east,font=\small] at (-0.05,2) {4};
\draw[thick] (-0.05,2.5) -- (0.05,2.5);
\node[anchor=east,font=\small] at (-0.05,2.5) {5};
\draw[thick] (-0.05,3) -- (0.05,3);
\node[anchor=east,font=\small] at (-0.05,3) {6};
\draw[thick] (-0.05,3.5) -- (0.05,3.5);
\node[anchor=east,font=\small] at (-0.05,3.5) {7};
\draw[thick] (-0.05,4) -- (0.05,4);
\node[anchor=east,font=\small] at (-0.05,4) {8};
\draw[thick] (0.5,-0.05) -- (0.5,0.05);
\node[anchor=north,font=\small] at (0.5,-0.05) {1};
\draw[thick] (1,-0.05) -- (1,0.05);
\node[anchor=north,font=\small] at (1,-0.05) {2};
\draw[thick] (1.5,-0.05) -- (1.5,0.05);
\node[anchor=north,font=\small] at (1.5,-0.05) {3};
\draw[thick] (2,-0.05) -- (2,0.05);
\node[anchor=north,font=\small] at (2,-0.05) {4};
\draw[thick] (2.5,-0.05) -- (2.5,0.05);
\node[anchor=north,font=\small] at (2.5,-0.05) {5};
\draw[thick] (3,-0.05) -- (3,0.05);
\node[anchor=north,font=\small] at (3,-0.05) {6};
\draw[thick] (3.5,-0.05) -- (3.5,0.05);
\node[anchor=north,font=\small] at (3.5,-0.05) {7};
\draw[thick] (4,-0.05) -- (4,0.05);
\node[anchor=north,font=\small] at (4,-0.05) {8};
\draw[very thick,red] (1,0) -- (4,4);
\filldraw [blue] (1,0) circle (2.0pt);
\filldraw [blue] (1.5,0) circle (2.0pt);
\filldraw [blue] (2,0) circle (2.0pt);
\filldraw [blue] (2.5,0.5) circle (2.0pt);
\filldraw [blue] (2.5,1) circle (2.0pt);
\filldraw [blue] (2.5,1.5) circle (2.0pt);
\filldraw [blue] (2.5,2) circle (2.0pt);
\filldraw [blue] (2,0.5) circle (2.0pt);
\filldraw [blue] (2,1) circle (2.0pt);
\filldraw [blue] (1.5,0.5) circle (2.0pt);
\filldraw [blue] (3,2.5) circle (2.0pt);
\filldraw [blue] (3.5,3) circle (2.0pt);
\filldraw [blue] (3.5,2.5) circle (2.0pt);
\filldraw [blue] (3.5,1.5) circle (2.0pt);
\filldraw [blue] (3.5,1) circle (2.0pt);
\filldraw [blue] (3,1) circle (2.0pt);
\filldraw [blue] (3,.5) circle (2.0pt);
\filldraw [blue] (3,0) circle (2.0pt);
\filldraw [blue] (2.5,0) circle (2.0pt);
\filldraw [blue] (3,1.5) circle (2.0pt);
\filldraw [blue] (3,2) circle (2.0pt);
\filldraw [blue] (3.5,2) circle (2.0pt);
\filldraw [blue] (4,3) circle (2.0pt);
\filldraw [blue] (4,4) circle (2.0pt);
\end{tikzpicture}
\caption{Subgraph plot of $\Sigma_3(C_4,0)$.} 
\end{subfigure}
\caption{Set $\Sigma_n(C_4,0)$ and upper boundary of its convex hull (in red) for $n=2,3$.}\label{fig2}
\end{figure}

\vspace{-0.3cm}

\noindent 
Figure~\ref{fig2} can be plotted via the following commands, after loading the SageMath code listed in \textcolor{black}{the} appendix.  

\begin{table}[H] 
\centering
\scriptsize 
 \resizebox{0.9\textwidth}{!}
    {
  \begin{tabular}{|ll|l|} 
 \hline
 \multicolumn{2}{|l|}{
   \EscVerb{G = [[1,2],[2,3],[3,4],[4,1]]; EP = []; SG2=convexhull(2,G,EP); SG3=convexhull(3,G,EP)
   }
 }
 \\
 \hline
 \multicolumn{2}{|l|}{
\EscVerb{Polyhedron(SG2).plot(color = "pink")+point(SG2,color = "blue",size=20)}
}
 \\
 \hline
 \multicolumn{2}{|l|}{
\EscVerb{Polyhedron(SG3).plot(color = "pink")+point(SG3,color = "blue",size=20)}
 }
 \\
\hline
\end{tabular}
}
\end{table} 
\end{example}
\vspace{-0.6cm}
\begin{example}
\noindent
 Let $G=C_4$ be a rooted $4$-cycle with
  one endpoint, i.e. $r=3$ and $m=1$,
  see Figure~\ref{fig5}.

\vspace{-0.6cm}

\begin{figure}[H]
  \centering
   \begin{tikzpicture}[scale=0.8] 
\draw[black, thick] (0,0) rectangle (3,3);
\filldraw [blue] (0,0) circle (3pt);
\filldraw [red] (0,3) circle (3pt);
\filldraw [blue] (3,0) circle (3pt);
\filldraw [blue] (3,3) circle (3pt);

\node[anchor=north,font=\small] at (0,0) {1};
\node[anchor=north,font=\small] at (3,0) {2};
\node[anchor=south,font=\small] at (3,3) {3};
\node[anchor=south,red,font=\small] at (0,3) {4};

   \end{tikzpicture}
   \caption{Connected graph $G=C_4$ with $r=3$ and $m=1$.}
   \label{fig5}
   \end{figure}

   \vspace{-0.3cm}

\noindent
 The sets 
  $\Sigma_3(C_4,1)$ and $\Sigma_4(C_4,1)$
  are plotted
  in Figure~\ref{fig6}. 

   \vspace{-0.2cm}

   \begin{figure}[H]
    \begin{subfigure}{.5\textwidth}
    \vskip2cm
   \centering
   \begin{tikzpicture}[scale=0.8] 
  \draw[step={(1cm,0.5cm)}, very thin, gray!40] (0,0) grid (7,4.5);
  \draw[->] (0,0) -- (0,4.5);
\node[anchor=east,font=\small] at (-0.05,4.5) {$y$};
\draw[->] (0,0) -- (7,0);
\node[anchor=north,font=\small] at (7,-0.05) {$x$};
\draw[thick] (-0.05,0.5) -- (0.05,0.5);
\node[anchor=east,font=\small] at (-0.05,0.5) {1};
\draw[thick] (-0.05,1) -- (0.05,1);
\node[anchor=east,font=\small] at (-0.05,1) {2};
\draw[thick] (-0.05,2) -- (0.05,2);
\node[anchor=east,font=\small] at (-0.05,2) {4};
\draw[thick] (-0.05,3) -- (0.05,3);
\node[anchor=east,font=\small] at (-0.05,3) {6};
\draw[thick] (-0.05,1.5) -- (0.05,1.5);
\node[anchor=east,font=\small] at (-0.05,1.5) {3};
\draw[thick] (-0.05,2.5) -- (0.05,2.5);
\node[anchor=east,font=\small] at (-0.05,2.5) {5};
\draw[thick] (-0.05,3) -- (0.05,3);
\node[anchor=east,font=\small] at (-0.05,3) {6};
\draw[thick] (-0.05,3.5) -- (0.05,3.5);
\node[anchor=east,font=\small] at (-0.05,3.5) {7};
\draw[thick] (-0.05,4) -- (0.05,4);
\node[anchor=east,font=\small] at (-0.05,4) {8};
\draw[thick] (1,-0.05) -- (1,0.05);
\node[anchor=north,font=\small] at (1,-0.05) {1};
\draw[thick] (2,-0.05) -- (2,0.05);
\node[anchor=north,font=\small] at (2,-0.05) {2};
\draw[thick] (3,-0.05) -- (3,0.05);
\node[anchor=north,font=\small] at (3,-0.05) {3};
\draw[thick] (4,-0.05) -- (4,0.05);
\node[anchor=north,font=\small] at (4,-0.05) {4};
\draw[thick] (5,-0.05) -- (5,0.05);
\node[anchor=north,font=\small] at (5,-0.05) {5};
\draw[thick] (6,-0.05) -- (6,0.05);
\node[anchor=north,font=\small] at (6,-0.05) {6};
\draw[very thick,red] (2,1) -- (6,4);
\filldraw [blue] (2,0) circle (2.0pt);
\filldraw [blue] (2,0.5) circle (2.0pt);
\filldraw [blue] (2,1) circle (2.0pt);
\filldraw [blue] (3,0) circle (2.0pt);
\filldraw [blue] (3,0.5) circle (2.0pt);
\filldraw [blue] (3,1) circle (2.0pt);
\filldraw [blue] (3,1.5) circle (2.0pt);
\filldraw [blue] (4,0.5) circle (2.0pt);
\filldraw [blue] (4,1) circle (2.0pt);
\filldraw [blue] (4,1.5) circle (2.0pt);
\filldraw [blue] (4,2) circle (2.0pt);
\filldraw [blue] (4,2.5) circle (2.0pt);
\filldraw [blue] (5,1) circle (2.0pt);
\filldraw [blue] (5,1.5) circle (2.0pt);
\filldraw [blue] (5,2) circle (2.0pt);
\filldraw [blue] (5,2.5) circle (2.0pt);
\filldraw [blue] (5,3) circle (2.0pt);
\filldraw [blue] (6,3) circle (2.0pt);
\filldraw [blue] (6,4) circle (2.0pt);
  \end{tikzpicture}
   \caption{Subgraph plot of $\Sigma_3(C_4,1)$.} 
  \end{subfigure} 
  \begin{subfigure}{.5\textwidth}
   \centering
 \begin{tikzpicture}[scale=0.8] 
    \draw[->] (0,0) -- (0,6.5);
\draw[step={(0.5cm,0.5cm)}, very thin, gray!40] (0,0) grid (5,6.5);
    \node[anchor=east,font=\small] at (-0.05,6.5) {$y$};
    \draw[->] (0,0) -- (5,0);
    \node[anchor=north,font=\small] at (5,-0.05) {$x$};
    \draw[thick] (-0.05,0.5) -- (0.05,0.5);
    \node[anchor=east,font=\small] at (-0.05,0.5) {1};
    \draw[thick] (-0.05,1) -- (0.05,1);
    \node[anchor=east,font=\small] at (-0.05,1) {2};
    \draw[thick] (-0.05,1.5) -- (0.05,1.5);
    \node[anchor=east,font=\small] at (-0.05,1.5) {3};
    \draw[thick] (-0.05,2) -- (0.05,2);
    \node[anchor=east,font=\small] at (-0.05,2) {4};
    \draw[thick] (-0.05,2.5) -- (0.05,2.5);
    \node[anchor=east,font=\small] at (-0.05,2.5) {5};
    \draw[thick] (-0.05,3) -- (0.05,3);
    \node[anchor=east,font=\small] at (-0.05,3) {6};
    \draw[thick] (-0.05,3.5) -- (0.05,3.5);
    \node[anchor=east,font=\small] at (-0.05,3.5) {7};
    \draw[thick] (-0.05,4) -- (0.05,4);
    \node[anchor=east,font=\small] at (-0.05,4) {8};
    \draw[thick] (-0.05,4.5) -- (0.05,4.5);
    \node[anchor=east,font=\small] at (-0.05,4.5) {9};
    \draw[thick] (-0.05,5) -- (0.05,5);
    \node[anchor=east,font=\small] at (-0.05,5) {10};
    \draw[thick] (-0.05,5.5) -- (0.05,5.5);
    \node[anchor=east,font=\small] at (-0.05,5.5) {11};
    \draw[thick] (-0.05,6) -- (0.05,6);
    \node[anchor=east,font=\small] at (-0.05,6) {12};
    \draw[thick] (0.5,-0.05) -- (0.5,0.05);
    \node[anchor=north,font=\small] at (0.5,-0.05) {1};
    \draw[thick] (1,-0.05) -- (1,0.05);
    \node[anchor=north,font=\small] at (1,-0.05) {2};
    \draw[thick] (1.5,-0.05) -- (1.5,0.05);
    \node[anchor=north,font=\small] at (1.5,-0.05) {3};
    \draw[thick] (2,-0.05) -- (2,0.05);
    \node[anchor=north,font=\small] at (2,-0.05) {4};
    \draw[thick] (2.5,-0.05) -- (2.5,0.05);
    \node[anchor=north,font=\small] at (2.5,-0.05) {5};
    \draw[thick] (3,-0.05) -- (3,0.05);
    \node[anchor=north,font=\small] at (3,-0.05) {6};
    \draw[thick] (3.5,-0.05) -- (3.5,0.05);
    \node[anchor=north,font=\small] at (3.5,-0.05) {7};
    \draw[thick] (4,-0.05) -- (4,0.05);
    \node[anchor=north,font=\small] at (4,-0.05) {8};
    \draw[thick] (4.5,-0.05) -- (4.5,0.05);
    \node[anchor=north,font=\small] at (4.5,-0.05) {9};
    \draw[very thick,red] (1.5,1.5) -- (4.5,6);
    \filldraw [blue] (1.5,0) circle (2.0pt);
    \filldraw [blue] (1.5,0.5) circle (2.0pt);
    \filldraw [blue] (1.5,1) circle (2.0pt);
    \filldraw [blue] (1.5,1.5) circle (2.0pt);
    \filldraw [blue] (2,0) circle (2.0pt);
    \filldraw [blue] (2,0.5) circle (2.0pt);
    \filldraw [blue] (2,1) circle (2.0pt);
    \filldraw [blue] (2,1.5) circle (2.0pt);
    \filldraw [blue] (2,2) circle (2.0pt);
    \filldraw [blue] (2.5,0.5) circle (2.0pt);
    \filldraw [blue] (2.5,1) circle (2.0pt);
    \filldraw [blue] (2.5,1.5) circle (2.0pt);
    \filldraw [blue] (2.5,2) circle (2.0pt);
    \filldraw [blue] (2.5,2.5) circle (2.0pt);
    \filldraw [blue] (2.5,3) circle (2.0pt);
    \filldraw [blue] (3,1) circle (2.0pt);
    \filldraw [blue] (3,1.5) circle (2.0pt);
    \filldraw [blue] (3,2) circle (2.0pt);
    \filldraw [blue] (3,2.5) circle (2.0pt);
    \filldraw [blue] (3,3) circle (2.0pt);
    \filldraw [blue] (3,3.5) circle (2.0pt);
    
    \filldraw [blue] (3.5,1.5) circle (2.0pt);
    \filldraw [blue] (3.5,2) circle (2.0pt);
    \filldraw [blue] (3.5,2.5) circle (2.0pt);
    \filldraw [blue] (3.5,3) circle (2.0pt);
    \filldraw [blue] (3.5,3.5) circle (2.0pt);
    \filldraw [blue] (3.5,4) circle (2.0pt);
    \filldraw [blue] (3.5,4.5) circle (2.0pt);

    \filldraw [blue] (4,3) circle (2.0pt);
    \filldraw [blue] (4,3.5) circle (2.0pt);
    \filldraw [blue] (4,4) circle (2.0pt);
    \filldraw [blue] (4,4.5) circle (2.0pt);
    \filldraw [blue] (4,5) circle (2.0pt);

    \filldraw [blue] (4.5,5) circle (2.0pt);
    \filldraw [blue] (4.5,6) circle (2.0pt);
  \end{tikzpicture}
 \caption{Subgraph plot of $\Sigma_4(C_4,1)$.} 
  \end{subfigure}
  \caption{Set $\Sigma_n(C_4,1)$ and upper boundary of its convex hull (in red) for $n=3,4$.}\label{fig6}
  \end{figure}
\vspace{-0.3cm}

\noindent 
 Figure~\ref{fig6} can be plotted via the following commands. 

\begin{table}[H] 
\centering
\scriptsize 
 \resizebox{0.9\textwidth}{!}
    {
  \begin{tabular}{|ll|l|} 
 \hline
 \multicolumn{2}{|l|}{
   \EscVerb{G = [[1,2],[2,3]]; EP = [[1,3]]; SG3=convexhull(3,G,EP); SG4=convexhull(4,G,EP)
   }
 }
 \\
 \hline
 \multicolumn{2}{|l|}{
\EscVerb{Polyhedron(SG3).plot(color = "pink")+point(SG3,color = "blue",size=20)}
}
 \\
 \hline
 \multicolumn{2}{|l|}{
\EscVerb{Polyhedron(SG4).plot(color = "pink")+point(SG4,color = "blue",size=20)}
 }
 \\
\hline
\end{tabular}
}
\end{table} 

\end{example}

\vspace{-0.2cm}

\noindent
 Figures~\ref{fig1}-\ref{fig6} also show 
 an upper boundary plotted in red, which is characterized
 in the next definition. 
 \begin{definition}\label{def:uppboundary}
 We let $\widehat{\Sigma}_n(G,m)$ denote the upper boundary of the convex hull of $\Sigma_n(G,m)$, with 
 \begin{equation}
   \label{jkfl000} 
   \widehat{\Sigma}_n(G,m)
   \cap
   {\Sigma}_n(G,m)=
   \{(x_i,z_i):i=0,1,\dots, l \}, 
\end{equation} 
 for some $l\geq 1$, where 
$$n-1 = x_0<x_1 < \cdots < x_l = (n-1)r < x_{l +1}:= + \infty. 
$$
\end{definition}
 \textcolor{black}{In Definition~\ref{xyplot}, $x ( \rho_G ) =nr+m-v(\rho_G)=nr-|\rho|$ stands for the number of vertices being removed  in the process of graph contraction. Later on, in Definition~\ref{def:uppboundary} $x_0$ and $x_l$
    will be used to denote the minimum and the maximum number of vertices being removed. Because for all $\rho\in {\rm CNF} (n,r)$, 
$$r\le|\rho|\le n(r-1)+1,$$
 we have $x_0=n-1$ and $x_l =(n-1)r.$}
 We note that for any point $(x_i,z_i)$ in
 $\widehat{\Sigma}_n(G,m)
   \cap
       {\Sigma}_n(G,m)$, $i\in \{0,1,\ldots , l\}$,
 there exists a connected non-flat partition $\rho_{i} \in {\rm CNF} (n,r)$
 such that the associated graph $\rho_{i,G}$ satisfies
 \begin{equation}
   \label{fjklf} 
v(\rho_{i,G})=nr+m-x_i \quad \mbox{and} \quad e(\rho_{i,G})=ne(G)-z_i. 
\end{equation}
 We note that the upper boundary $\widehat{\Sigma}_n(G,m)$
starts at
$(x_0,y_0):=(n-1,(n-1)a_m(G))$ and ends at
$((n-1)r,(n-1) e(G))$, 
where $a_m(G)$ is defined in \eqref{maxdegree}.

\medskip 

We also recall the following lemma from \cite[Lemma~2.8]{LiuPrivault},
in which the maximality of connected non-flat partitions refers to maximizing the number of blocks, see also Proposition~6.1 in \cite{schulte-thaele}. 
\begin{lemma}
\label{lma2}
  \noindent
  $a)$ The cardinality of the set\  ${\rm CNF} (n,r)$
 of connected non-flat partitions of $[n]\times[r]$ satisfies 
 \begin{equation}
   \label{coeff-0}
  |{\rm CNF} (n,r)| \leq n!^r r!^{n-1}, 
  \qquad n,r \geq 1. 
\end{equation}
\noindent
$b)$ 
Let $\mathcal{M}(n,r)$ denote the \textcolor{black}{set}
 of \textcolor{black}{maximal} connected non-flat partitions of $[n]\times[r]$.
 \textcolor{black}{Then, each element of $\mathcal{M}(n,r)$ has precisely $(n-1)r+1$ blocks, and we have}
 \begin{equation}
   \nonumber 
  |\mathcal{M}(n,r)|=r^{n-1}\prod_{i=1}^{n-1}(1+(r-1)i),
  \qquad n,r\geq 1, 
\end{equation}
 with the bounds 
\begin{equation}\label{coeff-2}
    ( (r-1)r )^{n-1}(n-1)!\le
    |\mathcal{M}(n,r)|
     \leq ( (r-1)r )^{n-1}n!, \quad n\geq 1, \ r\geq 2. 
\end{equation}
\end{lemma}
\section{Leading diagrams}
\label{leading}
\noindent 
\textcolor{black}{Based on the convex hull of
$\Sigma_n(G,m)$ given in Definition~\ref{xyplot}, in this section we identify the dominant asymptotic order and the leading contribution appearing in the expression \eqref{cum-eq} of $\kappa_n(N_G)$, $n\geq 1$, which is key to the derivation of normal approximation results via the cumulant method. }
\begin{definition}\label{defleading}
 \textcolor{black}{Given $G$ a connected graph} with $V_G=[r+m]$, a diagram $\rho \in {\rm CNF}(n,r)$, 
 $n \geq 1$, is said to be a leading diagram
 \textcolor{black}{for a given $(c_\lambda )_{\lambda >0}$}
 if, for every $\sigma \in {\rm CNF}(n,r)$ \textcolor{black}{satisfies that }
 \begin{equation}
   \nonumber 
 \lambda^{v(\sigma_G)}c_\lambda^{e(\sigma_G)}
 = O\big(\lambda^{v(\rho_G)}c_\lambda^{e(\rho_G)}\big),
  \qquad \mathrm{as}~\lambda\to\infty. 
\end{equation}
\end{definition}
  The characterization of leading diagrams will use
  the following definition. 
 \begin{definition}
 Let $\rho\in{\rm CNF}(n,r)$
 and
 $i_\rho \in \{0,\ldots , l \}$ \textcolor{black}{be} such that
     $$x_{i_\rho} \le x ( \rho_G ) <x_{i_\rho+1},
     $$
 where $x ( \rho_G )$ is given in Definition~\ref{xyplot}.
      Using the notation \eqref{jkfl000}, we define 
 $$
 \theta_- ( \rho_G ) : = \left\{
 \begin{array}{ll}
   +\infty, & 
   i_\rho=0,
   \medskip
   \\ 
   \displaystyle
   \frac{z_{i_\rho}-z_{i_\rho-1}}{x_{i_\rho}-x_{i_\rho-1}}, & 
   1 \leq i_\rho \leq l, 
 \end{array}
 \right.
 \quad \mbox{and} \quad 
 \theta_+ ( \rho_G )  : = 
 \left\{
 \begin{array}{ll}
   \displaystyle
   \frac{z_{i_\rho+1}-z_{i_\rho}}{x_{i_\rho+1}-x_{i_\rho}}, 
   & 
    0 \leq i_\rho < l , 
   \medskip
   \\ 
   0, & 
   i_\rho= l. 
 \end{array}
 \right.
$$ 
\end{definition} 
\noindent 
 We note the inequality 
 $$
 \theta_- ( \rho_G ) \geq \theta_+ ( \rho_G ),
 $$
 which holds because $\widehat{\Sigma}_n(G,m)$ is the
 upper boundary of the convex hull of $\Sigma_n(G,m)$. 
 We also say that a diagram $\rho\in {\rm CNF}(n,r)$ lies on the boundary $\widehat{\Sigma}_n(G,m)$ if the point $( x ( \rho_G ) , y ( \rho_G ) )$ does. 
Lemma~\ref{lm:slope}
states that leading diagrams can only lie on the upper boundary $\widehat{\Sigma}_n(G,m)$.
\begin{lemma}\label{lm:slope}
  Let $G$ be a \textcolor{black}{connected} graph with $V_G=[r+m]$, $r\geq 2$, $m\geq 0$. Let $n\geq 2$ and
  assume that $\lim_{\lambda \to \infty} c_\lambda =0$.
  \begin{enumerate}[1)]
  \item
    Every leading diagram $\rho \in {\rm CNF}(n,r)$ 
    lies on the
    upper boundary $\widehat{\Sigma}_n(G,m)$, i.e. 
  $$
  (nr+m-v(\rho_G),ne(G)-e(\rho_G)) \in \widehat{\Sigma}_n(G,m).
  $$
 \item
  If a diagram $\rho \in {\rm CNF}(n,r)$ lies on the upper boundary $\widehat{\Sigma}_n(G,m)$ and 
  \begin{equation}
\nonumber 
  \lambda c_\lambda^{ \theta_- ( \rho_G ) }=O(1) \ \mbox{ and } \
  \lambda c_\lambda^{ \theta_+ ( \rho_G ) }=\Omega(1), 
\end{equation}
 then $\rho$ is a leading diagram.
\end{enumerate}
\end{lemma}
\begin{Proof}
  \noindent
$1)$
Suppose that $\rho\in{\rm CNF}(n,r)$ does not lie on the boundary
$\widehat{\Sigma}_n(G,m)$, i.e. 
$$
( x ( \rho_G ) , y ( \rho_G ) ) :=
( nr+m-v(\rho_G) , ne(G)-e(\rho_G) )
\in\Sigma_n(G,m)\backslash \widehat{\Sigma}_n(G,m).
$$ 
 Using \eqref{jkfl000}, 
 if $x ( \rho_G ) =x_l$, we know that $y ( \rho_G ) <z_l,$ as $(x ( \rho_G ) ,y ( \rho_G ) )$ is not on the upper boundary $\widehat{\Sigma}_n(G,m)$. Therefore, we have 
 $$\lambda^{v(\rho_G)-v(\rho_{l,G})}c_\lambda^{e(\rho_G)-e (\rho_{l,G})}=c_\lambda^{z_{l}-y ( \rho_G )}\ll 1,
 $$
 and $\rho$ cannot be a leading diagram.

 \medskip

 If $x ( \rho_G ) <x_l$, we choose ${i_\rho}\in \{0,\ldots , l-1 \}$
such that
 $x_{i_\rho}\le x ( \rho_G ) <x_{{i_\rho}+1}$. 
Since $\widehat{\Sigma}_n(G,m)$ is the upper boundary of a
convex hull, by \eqref{fjklf} 
 if ${i_\rho}<l $ we have 
 $$\nu:=\frac{z_{{i_\rho}+1}-y ( \rho_G ) }{x_{{i_\rho}+1}-x ( \rho_G ) }
 > \theta_+ ( \rho_G ) = \frac{z_{{i_\rho}+1}-z_{i_\rho}}{x_{{i_\rho}+1}-x_{i_\rho}}
 = \frac{e(\rho_{{i_\rho},G})-e(\rho_{{i_\rho}+1,G})}{v(\rho_{{i_\rho},G})-v(\rho_{{i_\rho}+1,G})},
 $$
 i.e. $\lambda c_\lambda^\nu
 \ll\lambda c_\lambda^{\theta_+ ( \rho_G ) }$
 because $\lim_{\lambda \to \infty} c_\lambda =0$, and 
 we consider three cases.
 \begin{enumerate}[i)]
   \item If
 $
1\ll\lambda c_\lambda^\nu
$, we have 
\begin{align*}
  \lambda^{v({\rho_G})-v(\rho_{{i_\rho}+1,G})}c_\lambda^{e({\rho_G})-e(\rho_{{i_\rho}+1,G})}
  & =  \left(\lambda c_\lambda^\nu\right)^{v({\rho_G})-v(\rho_{{i_\rho}+1,G})}
  \\
  & \ll \left(\lambda c_\lambda^{\theta_+ (\rho_G )} \right)^{v({\rho_G})-v(\rho_{{i_\rho}+1,G})}
  \\
   &   \leq \left(\lambda c_\lambda^{\theta_+ (\rho_G )}\right)^{v(\rho_{{i_\rho},G})-v(\rho_{{i_\rho}+1,G})}
  \\
   & =  \lambda^{v(\rho_{{i_\rho},G})-v(\rho_{{i_\rho}+1,G})}c_\lambda^{e(\rho_{{i_\rho},G})-e(\rho_{{i_\rho}+1,G})}, 
\end{align*}
which implies $\lambda^{v({\rho_G})}c_\lambda^{e({\rho_G})}\ll \lambda^{v(\rho_{{i_\rho},G})}c_\lambda^{e(\rho_{{i_\rho},G})}$
 as $\lim_{\lambda \to \infty} c_\lambda =0$, hence
 the diagram $\rho$ is not leading. 
\item  If
 $\lim_{\lambda \to \infty} \lambda c_\lambda^\nu = c>0$, we have 
\begin{align*}
  \lambda^{v({\rho_G})-v(\rho_{{i_\rho}+1,G})}c_\lambda^{e({\rho_G})-e(\rho_{{i_\rho}+1,G})}
  & =  \left(\lambda c_\lambda^\nu\right)^{v({\rho_G})-v(\rho_{{i_\rho}+1,G})}
  \\
  & \ll \left(\lambda c_\lambda^{\theta_+ (\rho_G )} \right)^{v({\rho_G})-v(\rho_{{i_\rho}+1,G})}
  \\
   &  \asymp \left(\lambda c_\lambda^{\theta_+ (\rho_G )}\right)^{v(\rho_{{i_\rho},G})-v(\rho_{{i_\rho}+1,G})}
  \\
   & =  \lambda^{v(\rho_{{i_\rho},G})-v(\rho_{{i_\rho}+1,G})}c_\lambda^{e(\rho_{{i_\rho},G})-e(\rho_{{i_\rho}+1,G})}, 
\end{align*}
\noindent
 and we conclude as above. 
\item  If
 $\lambda c_\lambda^\nu\ll 1$, we have
 \begin{equation}
   \nonumber
   \lambda^{v({\rho_G})-v(\rho_{{i_\rho}+1,G})}c_\lambda^{e({\rho_G})-e(\rho_{{i_\rho}+1,G})}\ll 1, 
\end{equation}
 hence
 $$\lambda^{v({\rho_G})}c_\lambda^{e({\rho_G})}\ll \lambda^{v(\rho_{{i_\rho}+1,G})}c_\lambda^{e(\rho_{{i_\rho}+1,G})}.
 $$
\noindent 
 As a consequence of the above, we find  
\begin{equation}
\nonumber
\lambda^{v({\rho_G})}c_\lambda^{e({\rho_G})}
\ll
\lambda^{v(\rho_{{i_\rho},G})}c_\lambda^{e(\rho_{{i_\rho},G})}
  \ \mbox{ or } \ 
\lambda^{v({\rho_G})}c_\lambda^{e({\rho_G})}
\ll
\lambda^{v(\rho_{{i_\rho}+1,G})}c_\lambda^{e(\rho_{{i_\rho}+1,G})}, 
\end{equation}
 hence $\rho$ is not a leading diagram. 
\end{enumerate} 
\noindent
$2)$ Suppose that $\rho$ does lie on the boundary $\widehat{\Sigma}_{n}(G,m)$.
Then, there exists ${i_\rho}\in \{0,\dots,l\}$
 such that $(x ( \rho_{G} ) ,y ( \rho_{G} ) )=(x_{{i_\rho}},z_{{i_\rho}})$, and it holds that
$$
\lambda c_\lambda^{ \theta_- (\rho_G )}=O(1) \ \mbox{ and } \
  \lambda c_\lambda^{ \theta_+ (\rho_G ) }=\Omega(1). 
$$
  \begin{enumerate}[i)] 
  \item
    If $j<{i_\rho}$, then $x_{j}-x_{{i_\rho}}<0$ and
  $$
  \frac{z_{{i_\rho}}-z_{j}}{x_{{i_\rho}}-x_{j}}\ge\frac{z_{{i_\rho}}-z_{{i_\rho}-1}}{x_{{i_\rho}}-x_{{i_\rho}-1}}=\textcolor{black}{\theta_- (\rho_G)}.
  $$
  Hence, 
\begin{equation}
  \lambda c_{\lambda}^{(z_{{i_\rho}}-z_{j})/(x_{{i_\rho}}-x_{j})}=O\big(\lambda c_\lambda^{ \theta_- (\rho_G )}
  \big)=O(1).
\end{equation}
Now, since
\begin{eqnarray}
  \label{fjkl4}
  \frac{\lambda^{nr+m-x_{{i_\rho}}}c_{\lambda}^{ne(G)-y_{{i_\rho}}}}{\lambda^{nr+m-x_{j}}c_{\lambda}^{ne(G)-z_{j}}}=\lambda^{x_{j}-x_{{i_\rho}}}c_{\lambda}^{z_{j}-z_{{i_\rho}}}=\big(\lambda c_{\lambda}^{
    (z_{j}-z_{{i_\rho}})/(x_{j}-x_{{i_\rho}})}\big)^{x_{j}-x_{{i_\rho}}}, 
\end{eqnarray} 
 we find 
\begin{equation}
  \nonumber
  \frac{\lambda^{v(\rho_G)}c_{\lambda}^{e(\rho_G)}}{\lambda^{v(\rho_{j,G})}c_{\lambda}^{e(\rho_{j,G})}}=\Omega(1).
\end{equation}
\item
   If $j>{i_\rho}$, then $x_{j}-x_{{i_\rho}}>0$ and
$$
\frac{z_{j}-z_{{i_\rho}}}{x_{j}-x_{{i_\rho}}}\le \textcolor{black}{\theta_+(\rho_G )},
$$
 therefore
\begin{equation}
  \lambda c_{\lambda}^{ ( z_{{i_\rho}}-z_{j})/(x_{{i_\rho}}-x_{j})}=\Omega
 \big(\lambda c_\lambda^{ \theta_+ (\rho_G )}\big)=\Omega(1),
\end{equation}
and \eqref{fjkl4} shows that 
\begin{equation}
  \nonumber
  \frac{\lambda^{v(\rho_G)}c_{\lambda}^{e(\rho_G)}}{\lambda^{v(\rho_{j,G})}c_{\lambda}^{e(\rho_{j,G})}}=\Omega(1).
\end{equation}
  \end{enumerate}
  This ensures that $\rho$ is a leading diagram.
\end{Proof}
 Proposition~\ref{pp:mbal} provides a sufficient condition
 ensuring that the upper boundary $\widehat{\Sigma}_n(G,m)$ is a line segment, a property used in the proof of Corollary~\ref{themo-1-1}. 
\begin{prop}\label{pp:mbal}
  Let $G$ be a connected graph with $V_G=[r+m]$, $r\geq 2$, $m\geq 0$, such that 
  the balance condition~\eqref{mbalanced} holds. 
  Then, the upper boundary
  $\widehat{\Sigma}_n(G,m)$
  is a line segment for all $n\geq 1$. 
\end{prop}
\begin{Proof}
  To present the result in a more compact form, we will first show that the requirement that the upper boundary $\widehat{\Sigma}_n(G,m)$ is a line segment is equivalent to 
  \begin{equation}
    \label{mbalineq}
      \frac{e(G)-a_m(G)}{r-1}\leq \frac{e(\rho_G)-a_m(G)}{v(\rho_G)-m-1},
      \qquad
  \rho \in {\rm CNF} (n,r).
    \end{equation}
    Because the upper boundary $\widehat{\Sigma}_n(G,m)$ starts at $(x_0,z_0):=(n-1,(n-1)a_m(G))$ and ends at $(x_l,z_l):=((n-1)r,(n-1)e(G))$, the requirement that the upper boundary $\widehat{\Sigma}_n(G,m)$ is a line segment is equivalent to that for any $(x,z)\in\Sigma_n(G,m)$, 
    \begin{equation}
      \frac{z_l-z}{x_l-x}\ge\frac{z_l-z_0}{x_l-x_0}=\frac{e(G)-a_m(G)}{r-1}.
    \end{equation}
    Considering the Definition~\ref{xyplot}, we obtain that the requirement itself is further equivalent to for any diagram $\rho \in {\rm CNF}(n,r)$, 
    \begin{equation}
      \frac{e(\rho_G)-e(G)}{v(\rho_G)-m-r}\ge \frac{e(G)-a_m(G)}{r-1},
    \end{equation}
    and, after reorganizing, the above inequality
    becomes equivalent to \eqref{mbalineq}. It remains to show that the balance condition~\eqref{mbalanced} ensures that \eqref{mbalineq} is satisfied. Here, we apply an induction argument to see this. When $n=1$, the claim is trivial as the only element in ${\rm CNF} (1,r)$ is isomorphic to $G$. Suppose now that \eqref{mbalineq} holds up to the rank $n\geq 1$.
    Let $\rho\in{\rm CNF}(n+1,r)$ be a non-flat connected partition
    of $[n+1]\times[r]$ with associated graph $\rho_G$.
       By Lemma~\ref{restrict-partition}, up to reordering
  of $\{1,\ldots , n+1\}$ there exists
a partition
  $\sigma \in {\rm CNF}(n,r)$
 obtained by restriction of $\rho$
 to $[n]\times [r]$.
 Let $\sigma'$ denote the partition obtained by
 restriction of $\rho$ to $\{n+1\}\times[r]$,
 see Figure~\ref{diagram1}.
    
\vspace{-1.4cm}
  
\begin{figure}[H]
\captionsetup[subfigure]{font=footnotesize}
\centering
\subcaptionbox{Partition $\rho$ of $[n+1]\times [r]$.}[.3\textwidth]{
\begin{tikzpicture}[scale=0.9,hide labels]
\tikzstyle{VertexStyle}=[shape = circle, fill = blue!20, minimum size = 0pt, scale=0., text = white, hide labels]

\draw[black, thick] (0,0) rectangle (5,3.75);

\node[anchor=east,font=\small] at (0.8,3) {1};
\node[anchor=east,font=\small] at (0.8,2) {2};
\node[anchor=east,font=\small] at (0.8,1) {3};

\node[anchor=south,font=\small] at (1,0) {1};
\node[anchor=south,font=\small] at (2,0) {2};
\node[anchor=south,font=\small] at (3,0) {3};
\node[anchor=south,font=\small] at (4,0) {4};

\filldraw [black] (1,1) circle (2pt);
\filldraw [black] (2,1) circle (2pt);
\filldraw [black] (3,1) circle (2pt);
\filldraw [black] (4,1) circle (2pt);
\filldraw [black] (1,2) circle (2pt);
\filldraw [black] (2,2) circle (2pt);
\filldraw [black] (3,2) circle (2pt);
\filldraw [black] (4,2) circle (2pt);
\filldraw [black] (2,2) circle (2pt);
\filldraw [black] (1,3) circle (2pt);
\filldraw [black] (2,3) circle (2pt);
\filldraw [black] (3,3) circle (2pt);
\filldraw [black] (4,3) circle (2pt);

\draw[very thick] (1,2) -- (1,3); 
\draw[very thick] (2,1) -- (3,3); 

\node (1) [label=above:{}] at (1,5) {};
\node (2) [label=above:{}] at (2,5) {};
\node (3) [label=above:{}] at (3,5) {};
\node (4) [label=above:{}] at (4,5) {};

\node (5) [label=above:{}] at (1,4) {};
\node (6) [label=above:{}] at (2,4) {};
\node (7) [label=above:{}] at (3,4) {};
\node (8) [label=above:{}] at (4,4) {};

 \begin{pgfonlayer}{background}
    \filldraw [line width=4mm,black!3]
                            (0.2,0.25) rectangle (4.75,3.55);
  \end{pgfonlayer}
\end{tikzpicture}}
\hskip0.5cm
\subcaptionbox{Splitting of $[n+1]\times [r]$.}[.3\textwidth]{
\begin{tikzpicture}[scale=0.9] 
\draw[black, thick] (0,0) rectangle (5,3.75);

\node[anchor=east,font=\small] at (0.8,3) {1};
\node[anchor=east,font=\small] at (0.8,2) {2};
\node[anchor=east,font=\small] at (0.8,1) {3};

\node[anchor=south,font=\small] at (1,0) {1};
\node[anchor=south,font=\small] at (2,0) {2};
\node[anchor=south,font=\small] at (3,0) {3};
\node[anchor=south,font=\small] at (4,0) {4};

\filldraw [black] (1,1) circle (2pt);
\filldraw [black] (2,1) circle (2pt);
\filldraw [black] (3,1) circle (2pt);
\filldraw [black] (4,1) circle (2pt);
\filldraw [black] (1,2) circle (2pt);
\filldraw [black] (2,2) circle (2pt);
\filldraw [black] (3,2) circle (2pt);
\filldraw [black] (4,2) circle (2pt);
\filldraw [black] (2,2) circle (2pt);
\filldraw [black] (1,3) circle (2pt);
\filldraw [black] (2,3) circle (2pt);
\filldraw [black] (3,3) circle (2pt);
\filldraw [black] (4,3) circle (2pt);

\draw[very thick] (1,2) -- (1,3); 
\draw[very thick] (2,1) -- (3,3); 

\node (1) [label=above:{}] at (1,5) {};
\node (2) [label=above:{}] at (2,5) {};
\node (3) [label=above:{}] at (3,5) {};
\node (4) [label=above:{}] at (4,5) {};

\node (5) [label=above:{}] at (1,4) {};
\node (6) [label=above:{}] at (2,4) {};
\node (7) [label=above:{}] at (3,4) {};
\node (8) [label=above:{}] at (4,4) {};

\node (9) [label=above:{}] at (1,3) {};
\node (10) [label=above:{}] at (2,3) {};
\node (11) [label=above:{}] at (3,3) {};
\node (12) [label=above:{}] at (4,3) {};

\node (13) [label=above:{}] at (1,2) {};
\node (14) [label=above:{}] at (2,2) {};
\node (15) [label=above:{}] at (3,2) {};
\node (16) [label=above:{}] at (4,2) {};

\node (17) [label=above:{}] at (1,1) {};
\node (18) [label=above:{}] at (2,1) {};
\node (19) [label=above:{}] at (3,1) {};
\node (20) [label=above:{}] at (4,1) {};

\draw[very thick,purple] \convexpath{9,12,16,13}{.2cm};

\draw[very thick,purple] \convexpath{20,17,20,17}{.2cm};

 \begin{pgfonlayer}{background}
    \filldraw [line width=4mm,black!3]
                       (0.2,0.25) rectangle (4.75,3.55);
  \end{pgfonlayer}
\end{tikzpicture}}
\hskip0.5cm
\subcaptionbox{Partitions $\sigma$ and $\sigma'$.}[.3\textwidth]{
\begin{tikzpicture}[scale=0.9] 

\draw[black, thick] (0,0) rectangle (5,3.75);

\node[anchor=east,font=\small] at (0.8,3) {1};
\node[anchor=east,font=\small] at (0.8,2) {2};
\node[anchor=east,font=\small] at (0.8,1) {3};

\node[anchor=south,font=\small] at (1,0) {1};
\node[anchor=south,font=\small] at (2,0) {2};
\node[anchor=south,font=\small] at (3,0) {3};
\node[anchor=south,font=\small] at (4,0) {4};

\filldraw [black] (1,1) circle (2pt);
\filldraw [black] (2,1) circle (2pt);
\filldraw [black] (3,1) circle (2pt);
\filldraw [black] (4,1) circle (2pt);
\filldraw [black] (1,2) circle (2pt);
\filldraw [black] (2,2) circle (2pt);
\filldraw [black] (3,2) circle (2pt);
\filldraw [black] (4,2) circle (2pt);
\filldraw [black] (2,2) circle (2pt);
\filldraw [black] (1,3) circle (2pt);
\filldraw [black] (2,3) circle (2pt);
\filldraw [black] (3,3) circle (2pt);
\filldraw [black] (4,3) circle (2pt);

\draw[very thick] (1,2) -- (1,3);

\node (1) [label=above:{}] at (1,5) {};
\node (2) [label=above:{}] at (2,5) {};
\node (3) [label=above:{}] at (3,5) {};
\node (4) [label=above:{}] at (4,5) {};

\node (5) [label=above:{}] at (1,4) {};
\node (6) [label=above:{}] at (2,4) {};
\node (7) [label=above:{}] at (3,4) {};
\node (8) [label=above:{}] at (4,4) {};

\node (9) [label=above:{}] at (1,3) {};
\node (10) [label=above:{}] at (2,3) {};
\node (11) [label=above:{}] at (3,3) {};
\node (12) [label=above:{}] at (4,3) {};

\node (13) [label=above:{}] at (1,2) {};
\node (14) [label=above:{}] at (2,2) {};
\node (15) [label=above:{}] at (3,2) {};
\node (16) [label=above:{}] at (4,2) {};

\node (17) [label=above:{}] at (1,1) {};
\node (18) [label=above:{}] at (2,1) {};
\node (19) [label=above:{}] at (3,1) {};
\node (20) [label=above:{}] at (4,1) {};

\draw[very thick,purple] \convexpath{9,12,16,13}{.2cm};

\draw[very thick,purple] \convexpath{20,17,20,17}{.2cm};

 \begin{pgfonlayer}{background}
    \filldraw [line width=4mm,black!3]
                  (0.2,0.25) rectangle (4.75,3.55);
  \end{pgfonlayer}
\end{tikzpicture}}

\caption{Splitting of a partition $\rho$ into $\sigma$ and $\sigma'$ with $n=3$ and $r=4$.}
\label{diagram1}
\end{figure}

\vspace{-0.4cm}
  
\noindent
  Given a graph $G$ with $r$ vertices,
  let $\rho_G$ denote the graph 
  with vertex set $V(\rho_G)$ built on $\rho$
  as in Definition~\ref{defgraph-1},
  see Figure~\ref{diagram2} for an example
  with $G$ a graph on $r+m=6$ vertices including two endpoints\textcolor{black}{.}
    
\begin{figure}[H]
\captionsetup[subfigure]{font=footnotesize}
\centering
\subcaptionbox{Vertex set $V(\rho_G)$.}[.3\textwidth]{
\begin{tikzpicture}[scale=0.9] 
\draw[step=1cm, very thin, gray!40] (-0.75,0.25) grid (4.75,3.75);
\draw[black, thick] (-0.75,0.25) rectangle (4.75,3.75);

\filldraw [red] (-0.2,1.5) circle (2pt);
\filldraw [red] (-0.2,2.5) circle (2pt);

\filldraw [black] (1,1) circle (2pt);
\filldraw [black] (2,1) circle (2pt);
\filldraw [black] (3,1) circle (2pt);
\filldraw [black] (4,1) circle (2pt);
\filldraw [black] (1,2) circle (2pt);
\filldraw [black] (2,2) circle (2pt);
\filldraw [black] (3,2) circle (2pt);
\filldraw [black] (4,2) circle (2pt);
\filldraw [black] (2,2) circle (2pt);
\filldraw [black] (1,3) circle (2pt);
\filldraw [black] (2,3) circle (2pt);
\filldraw [black] (3,3) circle (2pt);
\filldraw [black] (4,3) circle (2pt);

\draw[very thick] (1,2) -- (1,3); 

\draw[very thick] (2,1) -- (3,3); 

\draw[thick,dash dot,blue] (-0.2,2.5) .. controls (1,3) .. (1,3);
\draw[thick,dash dot,blue] (-0.2,2.5) .. controls (1,1) .. (1,1);

\draw[thick,dash dot,blue] (-0.2,1.5) .. controls (1,3) .. (1,3);
\draw[thick,dash dot,blue] (-0.2,1.5) .. controls (1,1) .. (1,1);

\draw[thick,dash dot,blue] (-0.2,1.5) .. controls (1,2) .. (1,2);
\draw[thick,dash dot,blue] (-0.2,2.5) .. controls (1,2) .. (1,2);

\draw[thick,dash dot,blue] (1,3) .. controls (1.5,3) .. (2,3);
\draw[thick,dash dot,blue] (2,3) .. controls (2.5,3) .. (3,3);
\draw[thick,dash dot,blue] (3,3) .. controls (3.5,3) .. (4,3);

\draw[thick,dash dot,blue] (1,2) .. controls (1.5,2) .. (2,2);
\draw[thick,dash dot,blue] (2,2) .. controls (2.5,2) .. (3,2);
\draw[thick,dash dot,blue] (3,2) .. controls (3.5,2) .. (4,2);

\draw[thick,dash dot,blue] (1,1) .. controls (1.5,1) .. (2,1);
\draw[thick,dash dot,blue] (2,1) .. controls (2.5,1) .. (3,1);
\draw[thick,dash dot,blue] (3,1) .. controls (3.5,1) .. (4,1);

\draw[very thick,purple] \convexpath{9,12,16,13}{.2cm};

\draw[very thick,purple] \convexpath{20,17,20,17}{.2cm};

 \begin{pgfonlayer}{background}
    \filldraw [line width=4mm,black!3]
                     (-0.58,0.45) rectangle (4.55,3.55);
  \end{pgfonlayer}
\end{tikzpicture}}
\hskip2.5cm
\subcaptionbox{Graph $\rho_G$.}[.3\textwidth]{
\begin{tikzpicture}[scale=0.9] 
\draw[step=1cm, very thin, gray!40] (-0.75,0.25) grid (4.75,3.75);
\draw[black, thick] (-0.75,0.25) rectangle (4.75,3.75);

\filldraw [red] (-0.2,1.5) circle (2pt);
\filldraw [red] (-0.2,2.5) circle (2pt);

\filldraw [black] (1,1) circle (2pt);

\filldraw [black] (3,1) circle (2pt);
\filldraw [black] (4,1) circle (2pt);

\filldraw [black] (2,2) circle (2pt);
\filldraw [black] (3,2) circle (2pt);
\filldraw [black] (4,2) circle (2pt);
\filldraw [black] (2,2) circle (2pt);
\filldraw [black] (1,2.5) circle (2pt);
\filldraw [black] (2,3) circle (2pt);
\filldraw [black] (3,3) circle (2pt);
\filldraw [black] (4,3) circle (2pt);

\draw[thick,blue] (-0.2,2.5) .. controls (1,2.5) .. (1,2.5);
\draw[thick,blue] (-0.2,2.5) .. controls (1,1) .. (1,1);

\draw[thick,blue] (-0.2,1.5) .. controls (1,2.5) .. (1,2.5);
\draw[thick,blue] (-0.2,1.5) .. controls (1,1) .. (1,1);

\draw[thick,blue] (1,2.5) .. controls (2,3) .. (2,3);
\draw[thick,blue] (2,3) .. controls (2.5,3) .. (3,3);
\draw[thick,blue] (3,3) .. controls (3.5,3) .. (4,3);

\draw[thick,blue] (1,2.5) .. controls (2,2) .. (2,2);
\draw[thick,blue] (2,2) .. controls (2.5,2) .. (3,2);
\draw[thick,blue] (3,2) .. controls (3.5,2) .. (4,2);

\draw[thick,blue] (1,1) .. controls (1.5,1) .. (3,3);
\draw[thick,blue] (3,3) .. controls (2.5,1) .. (3,1);
\draw[thick,blue] (3,1) .. controls (3.5,1) .. (4,1);

\draw[very thick,purple] \convexpath{9,12,16,13}{.2cm};

\draw[very thick,purple] \convexpath{20,17,20,17}{.2cm};

 \begin{pgfonlayer}{background}
    \filldraw [line width=4mm,black!3]
                (-0.58,0.45) rectangle (4.55,3.55);
   \end{pgfonlayer}
\end{tikzpicture}}
\caption{Splitting of $\rho_G$ into $\sigma_G$ and $\sigma'_G$ with
  $n=3$, $r=4$ and two endpoints $m=2$.}
\label{diagram2}
\end{figure}

\vspace{-0.4cm}

\noindent
 Next, we consider the graphs $\sigma_G$ and $\sigma'_G$ 
 obtained from Definition~\ref{defgraph-1}
 on the vertex sets 
    $$
 V(\sigma_G):=
 \big\{
     b\in \rho \ : \ b \cap (\pi_1\cup \cdots \cup \pi_n) \not= \emptyset
     \big\} \cup [m]
     $$
     and
           $$
     V(\sigma'_G):=
     \big\{
     b\in \rho \ : \ b \cap \pi_{n+1} \not= \emptyset
     \big\} \cup [m], 
     $$
     with
     $\sigma'_G \simeq G$ because $\rho$ is non-flat,
     see Figure~\ref{diagram3}. 

     \vspace{-0.1cm}
  
\begin{figure}[H]
\captionsetup[subfigure]{font=footnotesize}
\centering
\subcaptionbox{Vertex set $V(\sigma_G)$.}[.3\textwidth]{
\begin{tikzpicture}[scale=0.9] 
\draw[step=1cm, very thin, gray!40] (-0.75,0.25) grid (4.75,3.75);
\draw[black, thick] (-0.75,0.25) rectangle (4.75,3.75);

\filldraw [red] (-0.2,1.5) circle (2pt);
\filldraw [red] (-0.2,2.5) circle (2pt);

\filldraw [black] (2,1) circle (2pt);

\filldraw [black] (1,2) circle (2pt);
\filldraw [black] (2,2) circle (2pt);
\filldraw [black] (3,2) circle (2pt);
\filldraw [black] (4,2) circle (2pt);
\filldraw [black] (2,2) circle (2pt);
\filldraw [black] (1,3) circle (2pt);
\filldraw [black] (2,3) circle (2pt);
\filldraw [black] (3,3) circle (2pt);
\filldraw [black] (4,3) circle (2pt);

\draw[very thick] (1,2) -- (1,3); 

\draw[very thick] (2,1) -- (3,3); 

\draw[thick,dash dot,blue] (-0.2,2.5) .. controls (1,3) .. (1,3);
\draw[thick,dash dot,blue] (-0.2,1.5) .. controls (1,3) .. (1,3);

\draw[thick,dash dot,blue] (-0.2,1.5) .. controls (1,2) .. (1,2);
\draw[thick,dash dot,blue] (-0.2,2.5) .. controls (1,2) .. (1,2);

\draw[thick,dash dot,blue] (1,3) .. controls (1.5,3) .. (2,3);
\draw[thick,dash dot,blue] (2,3) .. controls (2.5,3) .. (3,3);
\draw[thick,dash dot,blue] (3,3) .. controls (3.5,3) .. (4,3);

\draw[thick,dash dot,blue] (1,2) .. controls (1.5,2) .. (2,2);
\draw[thick,dash dot,blue] (2,2) .. controls (2.5,2) .. (3,2);
\draw[thick,dash dot,blue] (3,2) .. controls (3.5,2) .. (4,2);

\draw[very thick,purple] \convexpath{9,12,16,13}{.2cm};

 \begin{pgfonlayer}{background}
    \filldraw [line width=4mm,black!3]
              (-0.58,0.45) rectangle (4.55,3.55);
  \end{pgfonlayer}
\end{tikzpicture}}
\hskip0.5cm
\subcaptionbox{Vertex set $V(\sigma'_G)$.}[.3\textwidth]{
\begin{tikzpicture}[scale=0.9] 
\draw[step=1cm, very thin, gray!40] (-0.75,0.25) grid (4.75,3.75);
\draw[black, thick] (-0.75,0.25) rectangle (4.75,3.75);

\filldraw [red] (-0.2,1.5) circle (2pt);
\filldraw [red] (-0.2,2.5) circle (2pt);

\filldraw [black] (1,1) circle (2pt);
\filldraw [black] (2,1) circle (2pt);
\filldraw [black] (3,1) circle (2pt);
\filldraw [black] (4,1) circle (2pt);

\filldraw [black] (3,3) circle (2pt);

\draw[very thick] (2,1) -- (3,3); 

\draw[thick,dash dot,blue] (-0.2,2.5) .. controls (1,1) .. (1,1);

\draw[thick,dash dot,blue] (-0.2,1.5) .. controls (1,1) .. (1,1);

\draw[thick,dash dot,blue] (1,1) .. controls (1.5,1) .. (2,1);
\draw[thick,dash dot,blue] (2,1) .. controls (2.5,1) .. (3,1);
\draw[thick,dash dot,blue] (3,1) .. controls (3.5,1) .. (4,1);

\draw[very thick,purple] \convexpath{20,17,20,17}{.2cm};

 \begin{pgfonlayer}{background}
    \filldraw [line width=4mm,black!3]
       (-0.58,0.45) rectangle (4.55,3.55);
  \end{pgfonlayer}
\end{tikzpicture}}
\hskip0.5cm
\subcaptionbox{Vertex set $V(\sigma''_G)$.}[.3\textwidth]{
\begin{tikzpicture}[scale=0.9] 
\draw[step=1cm, very thin, gray!40] (-0.75,0.25) grid (4.75,3.75);
\draw[black, thick] (-0.75,0.25) rectangle (4.75,3.75);

\filldraw [red] (-0.2,1.5) circle (2pt);
\filldraw [red] (-0.2,2.5) circle (2pt);

\filldraw [black] (2,1) circle (2pt);

\filldraw [black] (3,3) circle (2pt);

\draw[very thick] (2,1) -- (3,3);

 \begin{pgfonlayer}{background}
    \filldraw [line width=4mm,black!3]
       (-0.58,0.45) rectangle (4.55,3.55);
  \end{pgfonlayer}
\end{tikzpicture}}
\hskip0.5cm
\caption{Splitting of $V(\rho_G)$ into $V(\sigma_G), V(\sigma'_G)$ with
  $n=3$, $r=4$ and two endpoints $m=2$.}
\label{diagram3}
\end{figure}

\vspace{-0.4cm}
     
\noindent
 Let now $\sigma''_G$ denote the graph induced by $\rho_G$ on 
   $V(\sigma''_G):= V(\sigma_G) \cap V(\sigma'_G)$, see
   Figure~\ref{diagram4}. 

\begin{figure}[H]
\captionsetup[subfigure]{font=footnotesize}
\centering
\subcaptionbox{Graph $\sigma_G$.}[.3\textwidth]{
\begin{tikzpicture}[scale=0.9] 
\draw[step=1cm, very thin, gray!40] (-0.75,0.25) grid (4.75,3.75);
\draw[black, thick] (-0.75,0.25) rectangle (4.75,3.75);

\filldraw [red] (-0.2,1.5) circle (2pt);
\filldraw [red] (-0.2,2.5) circle (2pt);

\filldraw [black] (1,2) circle (2pt);

\filldraw [black] (2,3) circle (2pt);
\filldraw [black] (3,3) circle (2pt);
\filldraw [black] (4,3) circle (2pt);

\filldraw [black] (2,1) circle (2pt);
\filldraw [black] (3,1) circle (2pt);
\filldraw [black] (4,1) circle (2pt);

\draw[thick,blue] (-0.2,1.5) .. controls (1,2) .. (1,2);
\draw[thick,blue] (-0.2,2.5) .. controls (1,2) .. (1,2);

\draw[thick,blue] (2,3) .. controls (2.5,3) .. (3,3);
\draw[thick,blue] (3,3) .. controls (3.5,3) .. (4,3);

\draw[thick,blue] (1,2) .. controls (2,1) .. (2,1);
\draw[thick,blue] (1,2) .. controls (2,3) .. (2,3);

\draw[thick,blue] (2,1) .. controls (2.5,1) .. (3,1);
\draw[thick,blue] (3,1) .. controls (3.5,1) .. (4,1);

 \begin{pgfonlayer}{background}
    \filldraw [line width=4mm,black!3]
        (-0.58,0.45) rectangle (4.55,3.55);
   \end{pgfonlayer}
\end{tikzpicture}}
\hskip0.5cm
\subcaptionbox{Graph $\sigma'_G$.}[.3\textwidth]{
\begin{tikzpicture}[scale=0.9] 
\draw[step=1cm, very thin, gray!40] (-0.75,0.25) grid (4.75,3.75);
\draw[black, thick] (-0.75,0.25) rectangle (4.75,3.75);

\filldraw [red] (-0.2,1.5) circle (2pt);
\filldraw [red] (-0.2,2.5) circle (2pt);

\filldraw [black] (1,2) circle (2pt);

\filldraw [black] (2,2) circle (2pt);
\filldraw [black] (3,2) circle (2pt);
\filldraw [black] (4,2) circle (2pt);

\draw[thick,blue] (-0.2,1.5) .. controls (1,2) .. (1,2);
\draw[thick,blue] (-0.2,2.5) .. controls (1,2) .. (1,2);

\draw[thick,blue] (2,2) .. controls (2.5,2) .. (3,2);
\draw[thick,blue] (3,2) .. controls (3.5,2) .. (4,2);

\draw[thick,blue] (1,2) .. controls (2,2) .. (2,2);

 \begin{pgfonlayer}{background}
    \filldraw [line width=4mm,black!3]
  (-0.58,0.45) rectangle (4.55,3.55);
 \end{pgfonlayer}
\end{tikzpicture}}
\hskip0.5cm
\subcaptionbox{Graph $\sigma''_G$.}[.3\textwidth]{
\begin{tikzpicture}[scale=0.9] 
\draw[step=1cm, very thin, gray!40] (-0.75,0.25) grid (4.75,3.75);
\draw[black, thick] (-0.75,0.25) rectangle (4.75,3.75);

\filldraw [red] (-0.2,1.5) circle (2pt);
\filldraw [red] (-0.2,2.5) circle (2pt);

\filldraw [black] (2.5,2) circle (2pt);

 \begin{pgfonlayer}{background}
    \filldraw [line width=4mm,black!3]
  (-0.58,0.45) rectangle (4.55,3.55);
 \end{pgfonlayer}
\end{tikzpicture}}
\caption{Splitting of $\rho_G$ into $\sigma_G$ and $\sigma'_G$ with
  $n=3$, $r=4$ and two endpoints $m=2$.}
\label{diagram4}
\end{figure}

\vspace{-0.4cm}

\noindent
  Then, $\sigma''_G$ contains $m$ endpoints
  in addition to at least one non-endpoint vertex
     due to the connectedness of $\rho$, hence
     we have $v(\sigma''_G)\ge m+1$. Since
     $\sigma''_G\subset \sigma'_G$ and $v(\sigma''_G)\ge m+1$,
     by the balance condition \eqref{mbalanced} we have 
       \begin{align*}
  \frac{e(\sigma''_G)-a_m(G)}{v(\sigma''_G)-m-1}\leq \frac{e(G)-a_m(G)}{r-1}, 
  \end{align*} 
       with the convention $0/0=0$.
       Hence, by the induction hypothesis \eqref{mbalineq} applied at the rank $n\geq 1$
       to $\sigma_G$, 
       we have 
 \begin{align*}
 & \frac{e(\rho_G)-a_m(G)}{v(\rho_G)-m-1}
   =
  \frac{(e(\sigma_G)-a_m(G))+(e(\sigma'_G)-a_m(G))-(e(\sigma''_G)-a_m(G))}{
    v(\sigma_G)+v(\sigma'_G)-v(\sigma''_G)-m-1}
  \\
  & \quad \quad
   =
  \frac{(e(\sigma_G)-a_m(G))+(e(G)-a_m(G))-(e(\sigma''_G)-a_m(G))}{v(\sigma_G)+v(G)-v(\sigma''_G)-m-1}
  \\
           & \quad \quad \geq 
  \frac{ ( v(\sigma_G)-m-1 ) \frac{e(G)-a_m(G)}{r-1}  
   +(e(G)-a_m(G))-(v(\sigma''_G)-m-1)\frac{e(G)-a_m(G)}{r-1}
   }{
    v(\sigma_G)+v(G)-v(\sigma''_G)-m-1}
  \\
           & \quad \quad = 
  \frac{e(G)-a_m(G)}{r-1}. 
  \end{align*} 
\end{Proof}
The balance condition~\eqref{mbalanced} 
turns out to be necessary in order to ensure the upper boundary
 $\widehat{\Sigma}_n(G,m)$ to be a line segment, as shown in the following counterexample.

\begin{counterexample}\label{counex}
\begin{enumerate}[]
\item 
  Consider the graph $G$ of Figure~\ref{fig:diagram5}, 
  which is not strongly balanced, with $r=4$ and $m=0$.   
\begin{figure}[H]

\centering
\begin{tikzpicture}
\filldraw [blue] (1,0) circle (2pt);
\filldraw [blue] (1,2) circle (2pt);
\filldraw [blue] (2.5,1) circle (2pt);
\filldraw [blue] (4,1) circle (2pt);
\draw[thick] (1,0) -- (1,2) -- (2.5,1) -- (1,0);
\draw[thick] (2.5,1) -- (4,1);
\end{tikzpicture}
\caption{Not strongly balanced graph $G$.}\label{fig:diagram5} 
\end{figure}

\vskip-0.5cm

  \noindent
 Figure~\ref{fig4} shows that the
  upper boundary 
  $\widehat{\Sigma}_n(G,m)$ is not a line segment
  for $n=2$ and $n=3$. 
  $\widehat{\Sigma}_2(G,0)$ and $\widehat{\Sigma}_3(G,0)$.
  
  \vspace{-0.4cm}
 
\begin{figure}[H]
\begin{subfigure}{.5\textwidth}
  \vskip1.5cm
   \centering
 \begin{tikzpicture}[scale=0.8] 
\draw[step={(1cm,0.5cm)}, very thin, gray!40] (0,0) grid (5,2.5);
\draw[->] (0,0) -- (0,2.5);
\node[anchor=east,font=\small] at (-0.05,2.5) {$y$};
\draw[->] (0,0) -- (5,0);
\node[anchor=north,font=\small] at (5,-0.05) {$x$};
\draw[thick] (-0.05,0.5) -- (0.05,0.5);
\node[anchor=east,font=\small] at (-0.05,0.5) {1};
\draw[thick] (-0.05,1) -- (0.05,1);
\node[anchor=east,font=\small] at (-0.05,1) {2};
\draw[thick] (-0.05,1.5) -- (0.05,1.5);
\node[anchor=east,font=\small] at (-0.05,1.5) {3};
\draw[thick] (-0.05,2) -- (0.05,2);
\node[anchor=east,font=\small] at (-0.05,2) {4};
\draw[thick] (1,-0.05) -- (1,0.05);
\node[anchor=north,font=\small] at (1,-0.05) {1};
\draw[thick] (2,-0.05) -- (2,0.05);
\node[anchor=north,font=\small] at (2,-0.05) {2};
\draw[thick] (3,-0.05) -- (3,0.05);
\node[anchor=north,font=\small] at (3,-0.05) {3};
\draw[thick] (4,-0.05) -- (4,0.05);
\node[anchor=north,font=\small] at (4,-0.05) {4};
\draw[very thick,gray] (1,0) -- (4,2);
\draw[very thick,red] (1,0) -- (3,1.5) -- (4,2);
\filldraw [blue] (4,1) circle (2.0pt);
\filldraw [blue] (3,0.5) circle (2.0pt);
\filldraw [blue] (3,0) circle (2.0pt);
\filldraw [blue] (3,1) circle (2.0pt);
\filldraw [blue] (2,0) circle (2.0pt);
\filldraw [blue] (2,0.5) circle (2.0pt);
\filldraw [blue] (1,0) circle (2.0pt);
\filldraw [blue] (4,1.5) circle (2.0pt);
\filldraw [blue] (4,2) circle (2.0pt);
\filldraw [blue] (3,1.5) circle (2.0pt);
\end{tikzpicture}
    \caption{Subgraph plot of $\Sigma_2(G,0)$ with $\widehat{\Sigma}_2(G,0)$.} 
\end{subfigure}
\begin{subfigure}{.5\textwidth}
  \centering
\begin{tikzpicture}[scale=0.8] 
\draw[step={(0.5cm,0.5cm)}, very thin, gray!40] (0,0) grid (4.5,4.5);
\draw[->] (0,0) -- (0,4.5);
\node[anchor=east,font=\small] at (-0.05,4.5) {$y$};
\draw[->] (0,0) -- (4.5,0);
\node[anchor=north,font=\small] at (4.5,-0.05) {$x$};
\draw[thick] (-0.05,0.5) -- (0.05,0.5);
\node[anchor=east,font=\small] at (-0.05,0.5) {1};
\draw[thick] (-0.05,1) -- (0.05,1);
\node[anchor=east,font=\small] at (-0.05,1) {2};
\draw[thick] (-0.05,1.5) -- (0.05,1.5);
\node[anchor=east,font=\small] at (-0.05,1.5) {3};
\draw[thick] (-0.05,2) -- (0.05,2);
\node[anchor=east,font=\small] at (-0.05,2) {4};
\draw[thick] (-0.05,2.5) -- (0.05,2.5);
\node[anchor=east,font=\small] at (-0.05,2.5) {5};
\draw[thick] (-0.05,3) -- (0.05,3);
\node[anchor=east,font=\small] at (-0.05,3) {6};
\draw[thick] (-0.05,3.5) -- (0.05,3.5);
\node[anchor=east,font=\small] at (-0.05,3.5) {7};
\draw[thick] (-0.05,4) -- (0.05,4);
\node[anchor=east,font=\small] at (-0.05,4) {8};
\draw[thick] (0.5,-0.05) -- (0.5,0.05);
\node[anchor=north,font=\small] at (0.5,-0.05) {1};
\draw[thick] (1,-0.05) -- (1,0.05);
\node[anchor=north,font=\small] at (1,-0.05) {2};
\draw[thick] (2,-0.05) -- (2,0.05);
\node[anchor=north,font=\small] at (2,-0.05) {4};
\draw[thick] (3,-0.05) -- (3,0.05);
\node[anchor=north,font=\small] at (3,-0.05) {6};
\draw[thick] (3.5,-0.05) -- (3.5,0.05);
\node[anchor=north,font=\small] at (3.5,-0.05) {7};
\draw[thick] (4,-0.05) -- (4,0.05);
\node[anchor=north,font=\small] at (4,-0.05) {8};
\draw[very thick,gray] (1,0) -- (4,4);
\draw[very thick,red] (1,0) -- (2,1.5) -- (3,3) -- (4,4);
\filldraw [blue] (4,3) circle (2.0pt);
\filldraw [blue] (3.5,1.5) circle (2.0pt);
\filldraw [blue] (3.5,2) circle (2.0pt);
\filldraw [blue] (3,1) circle (2.0pt);
\filldraw [blue] (3,1.5) circle (2.0pt);
\filldraw [blue] (2.5,1) circle (2.0pt);
\filldraw [blue] (3.5,2.5) circle (2.0pt);
\filldraw [blue] (2.5,0.5) circle (2.0pt);
\filldraw [blue] (3.5,1) circle (2.0pt);
\filldraw [blue] (3,0.5) circle (2.0pt);
\filldraw [blue] (2,0.5) circle (2.0pt);
\filldraw [blue] (3,2) circle (2.0pt);
\filldraw [blue] (2.5,0) circle (2.0pt);
\filldraw [blue] (3,0) circle (2.0pt);
\filldraw [blue] (2,0) circle (2.0pt);
\filldraw [blue] (3.5,3) circle (2.0pt);
\filldraw [blue] (2.5,1.5) circle (2.0pt);
\filldraw [blue] (2,1) circle (2.0pt);
\filldraw [blue] (3,2.5) circle (2.0pt);
\filldraw [blue] (1.5,0) circle (2.0pt);
\filldraw [blue] (1.5,0.5) circle (2.0pt);
\filldraw [blue] (2.5,2) circle (2.0pt);
\filldraw [blue] (1,0) circle (2.0pt);
\filldraw [blue] (2,1.5) circle (2.0pt);
\filldraw [blue] (4,3.5) circle (2.0pt);
\filldraw [blue] (4,4) circle (2.0pt);
\filldraw [blue] (3.5,3.5) circle (2.0pt);
\filldraw [blue] (3,3) circle (2.0pt);
\end{tikzpicture}
    \caption{Subgraph plot of $\Sigma_3(G,0)$ and $\widehat{\Sigma}_3(G,0)$.} 
\end{subfigure}
\caption{Set $\Sigma_n(G,0)$ and upper boundary of its convex hull (in red) for $n=2,3$.}
\label{fig4}
\end{figure}

\noindent 
 Figure~\ref{fig4} can be plotted via the following commands. 

\begin{table}[H] 
\centering
\scriptsize 
 \resizebox{0.9\textwidth}{!}
    {
  \begin{tabular}{|ll|l|} 
 \hline
 \multicolumn{2}{|l|}{
   \EscVerb{G = [[1,2],[2,3],[3,4],[1,3]]; EP = []; SG2=convexhull(2,G,EP); SG3=convexhull(3,G,EP)
   }
 }
 \\
 \hline
 \multicolumn{2}{|l|}{
\EscVerb{Polyhedron(SG2).plot(color = "pink")+point(SG2,color = "blue",size=20)}
}
 \\
 \hline
 \multicolumn{2}{|l|}{
\EscVerb{Polyhedron(SG3).plot(color = "pink")+point(SG3,color = "blue",size=20)}
 }
 \\
\hline
\end{tabular}
}
\end{table} 
\end{enumerate}

\vspace{-1.4cm}

\end{counterexample} 
 \textcolor{black}{
  \begin{remark}\label{counex2}
    We note that the balance condition \eqref{mbalanced}
    is not necessary for asymptotic normality of
    normalized subgraph counts.
 Consider the graph
 $G$ in Figure~\ref{fig8} for example, where $r=5, e(G)=7$, $m=a_m(G)=0$,
 which is not strongly balanced, and not even balanced. 
    \begin{figure}[H]
\begin{subfigure}{.5\textwidth}
\centering
\begin{tikzpicture}[scale=0.8]
\filldraw [blue] (1,0) circle (2pt);
\filldraw [blue] (1,2) circle (2pt);
\filldraw [blue] (2.5,1) circle (2pt);
\filldraw [blue] (-0.5,1) circle (2pt);
\filldraw [blue] (4,1) circle (2pt);
\draw[thick] (1,0) -- (1,2) -- (2.5,1) -- (1,0);
\draw[thick] (2.5,1) -- (4,1);
\draw[thick] (2.5,1) -- (-0.5,1);
\draw[thick] (-0.5,1) -- (1,2);
\draw[thick] (-.5,1) -- (1,0);
\end{tikzpicture}
\caption{ graph $G$.}\label{fig:diagram7} 
\end{subfigure}
\begin{subfigure}{.5\textwidth}
  \centering
\begin{tikzpicture}[scale=0.8]
  \filldraw [blue] (1,0) circle (2pt);
\filldraw [blue] (1,2) circle (2pt);
\filldraw [blue] (2.5,1) circle (2pt);
\filldraw [blue] (-.5,1) circle (2pt);
\draw[thick] (1,0) -- (1,2) -- (2.5,1) -- (1,0);
\draw[thick] (2.5,1) -- (-.5,1) -- (1,2);
\draw[thick] (-.5,1) -- (1,0);
\end{tikzpicture}
\caption{subgraph $H\subset G$.}\label{fig:diagram8} 
\end{subfigure}
\caption{A (not balanced) graph $G$ and subgraph $H$.}
\label{fig8}
\end{figure}
    \vspace{-.4cm}
    \noindent
 Since \eqref{mbalanced}
 is not satisfied,
 the upper boundary
 $\widehat{\Sigma}_n(G,m)$ is not a line segment,
 which leads to more potential candidates for leading diagram,
 beyond $\lambda^{1+(r-1)n}c_\lambda^{ne(G)-(n-1)a_m(G)}$ and $\lambda^r c_\lambda^{e(G)}$.
 Precisely, we have 
    \begin{eqnarray*}
      \kappa_n(N_G)       &\asymp&\max\left\{\lambda^{nr-(n-1)}c_{\lambda}^{ne(G)}, \lambda^{nr-(n-1)v(H)}c_\lambda^{ne(G)-(n-1)e(H)}, \lambda^{r} c_\lambda^{e(G)}\right\}\\
      &=&\max\left\{\lambda^{4n+1}c_\lambda^{7n}, \lambda^{n+4}c_{\lambda}^{n+6}, \lambda^5 c_\lambda^{7}\right\}\\
      &=& \left\{\begin{array}{ll}
    \lambda^{4n+1}c_\lambda^{7n} & \ \mbox{if }\ 
    c_\lambda
    \gg\lambda^{-1/2} ,
          \medskip
    \\    \displaystyle    \lambda^{n+4}c_{\lambda}^{n+6} & \ \mbox{if }\ \lambda^{-1}\ll c_\lambda\lesssim\lambda^{-1/2},     \medskip
      \\
      \displaystyle
      \lambda^{5}c_\lambda^{7} & \ \mbox{if }\  c_\lambda
      \lesssim \lambda^{-1}.
    \end{array}
    \right.
    \end{eqnarray*}
    Therefore, when $c_\lambda\gg\lambda^{-2/3}$
    we have $\kappa_n(\widebar{N}_G)\to0$, $n\ge3$, as $\lambda$ tends to infinity,
    which implies asymptotic normality of $\widebar{N}_G$ 
    by Theorem~1 in \cite{Janson1988}.
    On the other hand, when $c_\lambda\lesssim \lambda^{-1}$
    we have $\kappa_n(N_G)\to0$, $n\geq 1$,
    therefore $N_G$ does not have a Poisson limit.
  \end{remark}
 }
\section{Cumulant growth rates for subgraph counts} 
\label{cumulantrates}

\noindent
Under Assumptions~\ref{assm5-0}-\eqref{assm5-01-2}-\eqref{assm5-1-2}, 
    $F_\lambda (\rho)$ defined in
  \eqref{fjhkldf231} 
 satisfies 
 \begin{equation}
   \label{asymptics-1}
F_\lambda (\rho) \asymp \lambda^{|\rho|}c_\lambda^{e(\rho_G)}. 
\end{equation}
 In this section,
 we investigate the asymptotic behaviour of the cumulants
 $\kappa_n(N_G)$ in \eqref{cum-eq}
 as $c_\lambda\to0$ and the intensity $\lambda$ tends to infinity, 
 by identifying the leading diagram\textcolor{black}{s}
 $\rho\in{\rm CNF}(n,r)$
 which, from \eqref{cum-eq} and Definition~\ref{defleading},
 satisfy 
 \begin{equation}
   \label{eqvleading}
  \kappa_n(N_G) \asymp \lambda^{v(\rho_G)-m}c_\lambda^{e(\rho_G)}.
\end{equation}
  \begin{prop}\label{themo-1}
      Let $G$ be a connected graph with $V_G=[r+m]$ for $r\geq 2$ and $m$ endpoints,
    $m\geq 0$.
    Suppose that
    Assumptions~\ref{assm5-0}-\eqref{assm5-01-2}-\eqref{assm5-1-2}
    are satisfied
  and
    that the upper boundary $\widehat{\Sigma}_n(G,m)$ is a line segment
    linking $(n-1,(n-1)a_m(G))$ to $((n-1)r,(n-1)e(G))$. 
    \begin{enumerate}[a)]
      \item If $1 \gtrsim c_\lambda \gtrsim \lambda^{-( r-1 )/(e(G)-a_m(G)) }$, we have
    \begin{equation} \label{cas1}
      \kappa_n(N_G) \asymp \lambda^{1+(r-1)n}c_\lambda^{ne(G)-(n-1)a_m(G)},
    \quad \textcolor{black}{n \geq 2}.  
    \end{equation}
        \item If $c_\lambda \asymp \lambda^{-( r-1 )/(e(G)-a_m(G)) }$, we have
    \begin{equation} \label{cas3}
      \kappa_n(N_G)\asymp \lambda c_{\lambda}^{a_m(G)},
    \quad \textcolor{black}{n \geq 2}.   
  \end{equation}
  \item
    \label{partc}
    If $ c_\lambda \lesssim \lambda^{-( r-1 )/(e(G)-a_m(G)) }$, we find 
    \begin{equation} \label{cas2}
      \kappa_n(N_G) \asymp \lambda^r c_\lambda^{e(G)},
    \quad \textcolor{black}{n \geq 2}.  
    \end{equation}
    \end{enumerate} 
  \end{prop}
\begin{Proof}
 As in \eqref{jkfl000}, we write 
$$
 \widehat{\Sigma}_n(G,m) \cap
 {\Sigma}_n(G,m)
 =\{(x_0,z_0),(x_1,z_1),\dots,(x_l,z_l)\},
  $$
  with $x_0:=n-1<x_1<\cdots<x_l:=(n-1)r$.  
  According to Definition~\ref{xyplot}, we can find a corresponding
  partition $\rho_{i,G}\in{\rm CNF}(n,r)$ such that 
  $$
  \textcolor{black}{v(\rho_{i,G})=nr+m-x_i,~~~e(\rho_{i,G})=ne(G)-z_i.}
  $$
  Also, we write the (connected non-flat) set partition associated with $\rho_{i,G}$ as $\rho_i$,
  and from \eqref{asymptics-1} 
  we obtain that each $\rho_i$ contributes 
  \begin{equation}
    F_\lambda (\rho_i) \asymp \lambda^{v(\rho_{i,G})-m}c_\lambda^{e(\rho_{i,G})}
     \asymp \lambda^{nr-x_i}c_\lambda^{ne(G)-z_i}.
  \end{equation}
  Because the upper boundary is a line segment with endpoints $(n-1,(n-1)a_m(G))$ and $((n-1)r,(n-1)e(G))$, the slope of this line segment is
  $\theta:=(e(G)-a_m(G))/(r-1)$. By Lemma~\ref{lm:slope}-$(1)$
 the leading diagram $\rho$ must lie on 
 $\widehat{\Sigma}_n(G,m)$.  
 \begin{enumerate}[a)]
   \item 
     For any $j=1,\dots,l$, by
     \eqref{fjklf} 
 we have
\begin{eqnarray*} 
 \frac{\lambda^{v(\rho_{0,G})-m}c_\lambda^{e(\rho_{0,G})}}{\lambda^{v(\rho_{j,G})-m}c_\lambda^{e(\rho_{j,G})}} & = & 
 \frac{\lambda^{1+(r-1)n}c_\lambda^{ne(G)-(n-1)a_m(G)}}{\lambda^{nr-x_j}c_\lambda^{ne(G)-z_j}}
 \\
  & = & 
 \frac{\lambda^{nr-x_0}c_\lambda^{ne(G)-z_0}}{\lambda^{nr-x_j}c_\lambda^{ne(G)-z_j}}
 \\
  & = & 
 \lambda^{x_j-x_0} c_\lambda^{z_j-z_0}
 \\
  & = & 
 \lambda^{x_j-x_0} c_\lambda^{\theta (x_j-x_0)}
 \\
  & = & 
\big(\lambda c_\lambda^{(e(G)-a_m(G)) / ( r-1 )}\big)^{x_j-x_0},
\end{eqnarray*} 
    hence if $\lambda c_\lambda^{(e(G)-a_m(G))/(r-1)} \textcolor{black}{\gtrsim} 1$,
    we find 
$$
 \frac{\lambda^{v(\rho_{0,G})-m}c_\lambda^{e(\rho_{0,G})}}{\lambda^{v(\rho_{j,G})-m}c_\lambda^{e(\rho_{j,G})}}\textcolor{black}{\gtrsim} 1, 
  $$
 therefore any $\rho_G$ such that 
 $$
 (v(\rho_G)\textcolor{black}{-m},e(\rho_G))=(v(\rho_{0,G})-m,e(\rho_{0,G}))=(1+(r-1)n,ne(G)-(n-1)a_m(G))$$
 is a leading diagram, 
 and this yields \eqref{cas1} by \eqref{eqvleading}.  
 \item
 If $\lambda c_\lambda^{(e(G)-a_m(G))/(r-1)}\asymp 1$
  then any diagram $\rho_i$, $i=0,1, \ldots , l$
 on the segment $\widehat{\Sigma}_n(G,m)$ 
 is a leading diagram by Lemma~\ref{lm:slope}-$(2)$. 
 Furthermore, by choosing $j=l$ with
 $$
 (v(\rho_G)-m,e(\rho_G))=(v(\rho_{l , G})-m,e(\rho_{l , G}))=(r,e(G)),
 $$ 
 \textcolor{black}{we find that \eqref{eqvleading} yields \eqref{cas3}, i.e.} 
$$
 \kappa_n(N_G) \asymp \lambda^r c_\lambda^{e(G)}
  \asymp \textcolor{black}{\lambda c_{\lambda}^{a_m(G)}.}
$$ 
\item
 For any $j=0,\dots,l-1$, by \eqref{fjklf} 
     we have
\begin{eqnarray*} 
  \frac{\lambda^{v(\rho_{l , G})-m}c_\lambda^{e(\rho_{l , G})}}{\lambda^{v(\rho_{j , G})-m}c_\lambda^{e(\rho_{j , G})}}
  & = &
  \frac{\lambda^r c_\lambda^{e(G)}}{\lambda^{nr-x_j}c_\lambda^{ne(G)-z_j}}
  \\
  & = &
  \frac{\lambda^{nr-x_l}c_\lambda^{ne(G)-z_l}}{\lambda^{nr-x_j}c_\lambda^{ne(G)-z_j}}
  \\
  & = &
  \lambda^{x_j-x_l}c_\lambda^{z_j-z_l}
  \\
  & = &
  \lambda^{x_j-x_l}c_\lambda^{\theta (x_j-x_l )}
  \\
  & = &
  \big(\lambda c_\lambda^{(e(G)-a_m(G)) / ( r-1) }\big)^{x_j-x_l}, 
  \end{eqnarray*}
 hence if $\lambda c_\lambda^{(e(G)-a_m(G))/(r-1)} \textcolor{black}{\lesssim} 1$,
    we find 
$$
     \frac{\lambda^{v(\rho_{l , G})-m}c_\lambda^{e(\rho_{l , G})}}{\lambda^{v(\rho_{j , G})-m}c_\lambda^{e(\rho_{j , G})}}\textcolor{black}{\gtrsim} 1. 
  $$
 Therefore, any $\rho_G$ such that 
 $$
 (v(\rho_G)-m,e(\rho_G))=(v(\rho_{l , G})-m,e(\rho_{l , G}))=(r,e(G))
 $$
 is a leading diagram, 
 and this yields \eqref{cas2} by \eqref{eqvleading}.  

  \vspace{-0.5cm}

\end{enumerate} 
\end{Proof}
\noindent
 We note from Lemma~\ref{lm:slope} and Proposition~\ref{pp:mbal} 
 that as long as a connected graph $G$ satisfies the balance
 condition~\eqref{mbalanced}, the leading asymptotic order in
the expression \eqref{cum-eq}
of $\kappa_n(N_G)$ is fully determined by either the maximal or the minimal
connected non-flat partition.
Here, maximality, resp. minimality,
of partitions refers to the maximality,
resp. minimality, of their block counts. 
 As a consequence of Remark~\ref{jlkdfed23}, we have the following. 
\begin{remark}
  When $m=1$ with $a_1(G)=1$,
  Proposition~\ref{themo-1} holds for trees, cycles and complete graphs
  as they are all $K_2$-balanced 
  and the balance condition~\eqref{mbalanced} coincides with
    the $K_2$-balance condition \eqref{k2baldef}.
\end{remark}
 \noindent 
 By Propositions~\ref{pp:mbal} and \ref{themo-1},
 we have the following result. 
\begin{corollary}
\label{themo-1-1}
Let $G$ be a connected graph with $V_G=[r+m]$,
 $r\geq 2$, and $m$ endpoints,
$m\geq 0$. Suppose that
Assumptions~\ref{assm5-0}-\eqref{assm5-01-2}-\eqref{assm5-1-2}
  and the balance condition~\eqref{mbalanced} are satisfied.  
 Then, the cumulant $\kappa_n(N_G)$ of order $n\geq 1$ of
 the subgraph count $N_G$ satisfies the following. 
\begin{enumerate}[a)]
\item If $1 \gtrsim c_\lambda \gtrsim \lambda^{-( r-1 )/(e(G)-a_m(G)) }$,
  we have
    \begin{equation} \nonumber 
      \kappa_n(N_G) \asymp \lambda^{1+(r-1)n}c_\lambda^{ne(G)-(n-1)a_m(G)}. 
    \end{equation}
        \item If $c_\lambda \asymp \lambda^{-( r-1 )/(e(G)-a_m(G)) }$, we have
    \begin{equation} \nonumber 
      \kappa_n(N_G)\asymp \lambda c_{\lambda}^{a_m(G)}. 
  \end{equation}
  \item
    If $ c_\lambda \lesssim \lambda^{-( r-1 )/(e(G)-a_m(G)) }$, we find 
    \begin{equation}
      \nonumber 
      \kappa_n(N_G) \asymp \lambda^r c_\lambda^{e(G)}.
    \end{equation}
    \end{enumerate} 
  \end{corollary} 
 In addition, by Remark~\ref{jlkdfed23}-$a)$
   we have the following consequence of
  Corollary~\ref{themo-1-1}. 
  \begin{corollary}
     Let $G$ be a strongly balanced connected
 graph with $v(G)=r$ vertices, $r \geq 2$,
 and no endpoints, i.e. $m=0$,
 and suppose that
 Assumptions~\ref{assm5-0}-\eqref{assm5-01-2}-\eqref{assm5-1-2}
 are satisfied. 
\begin{enumerate}[a)]
  \item If $1 \gtrsim c_\lambda \gtrsim \lambda^{-( v(G)-1 )/e(G)} $, we have
\begin{equation}
\nonumber 
  \kappa_n(N_G) \asymp \lambda^{1+(v(G)-1)n}c_\lambda^{ne(G)}. 
\end{equation}
\item If $c_\lambda \asymp \lambda^{-( v(G)-1 )/e(G)}$, we have
\begin{equation}
\nonumber 
  \kappa_n(N_G)\asymp \lambda.  
\end{equation}
\item If $c_\lambda \lesssim \lambda^{-( v(G)-1 )/e(G)}$, we find 
\begin{equation}
\nonumber 
  \kappa_n(N_G)\asymp \lambda^{v(G)} c_\lambda^{e(G)}.
\end{equation}
\end{enumerate} 
  \end{corollary}
Proposition~\ref{thmstatu0}
deals with the cumulant growth of normalized subgraph counts,
for use in normal approximation. 
\textcolor{black}{In the particular case 
 of $(r+1)$-hop counting with $c_\lambda = 1$
 and $m=2$, \eqref{dkddd1} is consistent with the
 normalized cumulant bound (8.2) in \cite{privaultkhops}.} 
\begin{prop} 
\label{thmstatu0}
 Let $G$ be a connected graph with $V_G=[r+m]$ for $r\geq 2$ and $m\geq 0$.
   Suppose that
    Assumptions~\ref{assm5-0}-\eqref{assm5-01-2}-\eqref{assm5}
    are satisfied and
    that the upper boundary $\widehat{\Sigma}_n(G,m)$ is a line segment
    linking $(n-1,(n-1)a_m(G))$ to $((n-1)r,(n-1)e(G))$. 
  Denoting by
  $$
  \widebar{N}_G:=\frac{N_G-\kappa_1(N_G)}{\sqrt{\kappa_2(N_G)}}$$
  the \textcolor{black}{normalized} subgraph count, we have
  \begin{equation}
    \label{dkddd1}
     |\kappa_n(\widebar{N}_G)|\leq \frac{n!^r }{\Delta_\lambda^{n-2}},
     \quad \textcolor{black}{n \geq 3}, 
  \end{equation}
  where 
  \begin{eqnarray}
    \Delta_\lambda\asymp\left\{\begin{array}{ll}
    \lambda^{1/2}c_\lambda^{a_m(G)/2} & \ \mbox{if }\  \textcolor{black}{1\gtrsim}
    c_\lambda
    \gtrsim
    \lambda^{-(r-1)/(e(G)-a_m(G))} ,
          \medskip
    \\    \displaystyle    \lambda^{(e(G)-ra_m(G))/(2(e(G)-a_m(G)))} & \ \mbox{if }\  c_\lambda\asymp\lambda^{-(r-1)/(e(G)-a_m(G))},     \medskip
      \\
      \displaystyle
      \lambda^{r/2}c_\lambda^{e(G)/2} & \ \mbox{if }\ \lambda^{-r/e(G)}\ll c_\lambda
      \lesssim \lambda^{-(r-1)/(e(G)-a_m(G))}.
    \end{array}
    \right.
  \end{eqnarray}
\end{prop}
\begin{Proof}
  We start by assuming that $m\geq 1$. \textcolor{black}{We only focus on the case when $n\ge3$, as cases $n=1,2$ are trivial.}
  \begin{enumerate}[a)]
  \item When $\textcolor{black}{1\gtrsim}c_\lambda\textcolor{black}{\gtrsim}\lambda^{-(r-1)/(e(G)-a_m(G))}$, since
    the balance condition \eqref{mbalanced} holds, from
    Proposition~\ref{themo-1}-a), we know that  the leading diagrams belong to
\textcolor{black}{the set $\mathcal{M}(n,r)$ of 
  maximal connected non-flat partitions of $[n]\times[r]$,
   see Lemma~\ref{lma2}}. Therefore, given \eqref{cum-eq} and \eqref{coeff-0}, we can bound $\kappa_n(N_G)$ from above 
    \begin{eqnarray}\label{cumb-1}
      \kappa_n(N_G)&\leq  |{\rm CNF}(n,r)|\lambda^{1+(r-1)n}c_{\lambda}^{ne(G)-(n-1)a_m(G)} C_{1,n}\nonumber\\
      &\leq  n!^r r!^{n-1}\lambda^{1+(r-1)n}c_{\lambda}^{ne(G)-(n-1)a_m(G)} C_{1,n}
    \end{eqnarray}
    where
    \begin{equation}
      \label{fjkls}
      C_{1,n}:=\max_{\rho\in\mathcal{M}(n,r)}\int_{(\R^d)^{|\rho|}}\prod_{\substack{ 
      1 \leq j \leq m
      \\ i\in {\cal A}^\rho_j}}
    \varphi(x_i,y_j)
    \ \prod_{
      \substack{1 \leq k \textcolor{black}{<} l \le|\rho|
        \\
        \{ k , l \}\in E(\rho_G) 
    }}\varphi(x_k,x_l) \ \! \mathrm{d}x_1 \cdots\mathrm{d}x_{|\rho|}.
    \end{equation}
    Because the function $\varphi:\R^d\times\R^d\to[0,1]$ is symmetric and translation invariant, we can further bound $C_{1,n}$ as follows. From Definition~\ref{defgraph-1}, we know that for any $\rho\in\mathcal{M}(n,r)$, $\rho_G$ is a connected graph with $V(\rho_G)=[m+1+n(r-1)]$, as $|\rho|=1+n(r-1)$. Let $\tilde{\rho}_G$ be the subgraph of $\rho_G$ induced by $V(\tilde{\rho}_G)=[n(r-1)+2]$.
    And we also denote $\widebar{\rho}_G$ a spanning tree of $\tilde{\rho}_G$. Therefore, 
    \begin{align}
      \nonumber
      &\int_{(\R^d)^{|\rho|}}\prod_{\substack{ 
      1 \leq j \leq m
      \\ i\in {\cal A}^\rho_j}}
    \varphi(x_i,y_j)
    \ \prod_{
      \substack{1 \leq k \textcolor{black}{<} l \le|\rho|
        \\
        \{ k , l \}\in E(\rho_G) 
    }}\varphi(x_k,x_l) \ \! \mathrm{d}x_1 \cdots\mathrm{d}x_{|\rho|}
    \\
    \nonumber
    & \quad \leq \int_{(\R^d)^{|\rho|}}\prod_{
    i\in {\cal A}^\rho_1}
  \varphi(x_i,y_1)
  \ \prod_{
    \substack{1 \leq k \textcolor{black}{<} l \le|\rho|
      \\
      \{ k , l \}\in E(\tilde{\rho}_G) 
  }}\varphi(x_k,x_l) \ \! \mathrm{d}x_1 \cdots\mathrm{d}x_{|\rho|}
  \\
  \nonumber
  & \quad \le \int_{(\R^d)^{|\rho|}}
  \ \prod_{
    \substack{1 \leq k \textcolor{black}{<} l \le|\rho|+1
      \\
      \{ k , l \}\in E(\widebar{\rho}_G) 
  }}\varphi(x_k,x_l) \ \! \mathrm{d}x_1 \cdots\mathrm{d}x_{|\rho|}
  \\
  \label{fjkldf2}
   & \quad =\left(\int_{\R^d}\varphi(0,x) \ \! \mathrm{d}x\right)^{1+n(r-1)},
    \end{align} 
    where $x_{|\rho|+1}:=y_1$, and the last equality is obtained by
    integrating successively on the variables which
 correspond to leaves of $\widebar{\rho}_G$ as in the proofs of e.g.
 Theorem~7.1 of \cite{LNS21} or
 Lemma~3.1 of \cite{can2022} 
 since $\varphi$ is translation invariant
 by Assumption~\ref{assm5-0}-\eqref{assm5}.
 Hence, $C_{1,n}$ is bounded by $\zeta_\varphi^{1+n(r-1)}$, where 
 $$\zeta_\varphi:=\max
 \left(
 1 , \int_{\R^d}\varphi(0,x) \ \! \mathrm{d}x \right).
    $$ 
    In the other direction, as in \eqref{cum-eq}, \textcolor{black}{the} cumulants of $N_G$ \textcolor{black}{are} written as a summation of some non-negative terms.
    Therefore, we can bound $\kappa_n(N_G)$ from below, as follows: 
\begin{eqnarray}\label{cumb-2}
  \kappa_n(N_G)&\ge& |\mathcal{M}(n,r)|\lambda^{1+(r-1)n}c_{\lambda}^{ne(G)-(n-1)a_m(G)} C_{2,n}\nonumber\\
&\ge& ( (r-1)r )^{n-1}(n-1)!\lambda^{1+(r-1)n}c_{\lambda}^{ne(G)-(n-1)a_m(G)} C_{2,n},
\end{eqnarray}
where the last inequality comes from \eqref{coeff-2}, and 
\begin{equation}
  C_{2,n}:=\min_{\rho\in\mathcal{M}(n,r)}\int_{(\R^d)^{|\rho|}}\prod_{\substack{ 
      1 \leq j \leq m
      \\ i\in {\cal A}^\rho_j}}
    \varphi(x_i,y_j)
    \ \prod_{
      \substack{1 \leq k \textcolor{black}{<} l \le|\rho|
        \\
        \{ k , l \}\in E(\rho_G) 
    }}\varphi(x_l, x_k) \ \! \mathrm{d}x_1 \cdots\mathrm{d}x_{|\rho|}.
\end{equation}
Combining \eqref{cumb-1} and \eqref{cumb-2}, we have, for $n\ge3$
\begin{eqnarray}
 \kappa_n(\widebar{N}_G)&=& \frac{\kappa_n(N_G)}{\kappa_2(N_G)^{n/2}}\nonumber\\
  &\leq & \frac{n!^r r!^{n-1}\lambda^{1+(r-1)n}c_{\lambda}^{ne(G)-(n-1)a_m(G)} C_{1,n}}{
    \big( (r-1)r \lambda^{1+2(r-1)}c_{\lambda}^{2e(G)-a_m(G)} C_{2,2}\big)^{n/2}}\nonumber\\
  &=&n!^r\frac{r!^{n-1}}{((r-1)r)^{n/2}}\big(\lambda c_\lambda^{a_m(G)}\big)^{-(n-2)/2} \frac{C_{1,n}}{C_{2,2}^{n/2}}\nonumber\\
  &\leq & n!^r ((r-2)!)^{n-1}((r-1)r)^{n/2-1}\big( \lambda c_\lambda^{a_m(G)}\big)^{-(n-2)/2} \frac{\zeta_\varphi^{1+n(r-1)}}{C_{2,2}^{n/2}}\nonumber\\
  &\leq & \frac{n!^r}{\Big(C_3\sqrt{\lambda c_\lambda^{a_m(G)}}\Big)^{n-2}}\ ,
\end{eqnarray}
where $C_3$ is a constant depending on $r$ and $\varphi$.
\item When $c_\lambda\lesssim \lambda^{-(r-1)/(e(G)-a_m(G))}$,
  from Proposition~\ref{themo-1}-$b)$, the leading diagrams are $\rho\in{\rm CNF}(n,r)$ such that $\rho_G\simeq G$, which allows us to bound $\kappa_n(N_G)$ as follows:
    \begin{eqnarray}
      \kappa_n(N_G)&\leq & |{\rm CNF}(n,r)|\lambda^r c_{\lambda}^{e(G)} \int_{(\R^d)^r }
    \ \prod_{
      \substack{1 \leq k \textcolor{black}{<} l \le r+m
        \\
        \{ k , l \}\in E_G 
    }}\varphi(x_k,x_l) \ \! \mathrm{d}x_1 \cdots\mathrm{d}x_r\nonumber\\
    &\leq & n!^rr!^{n-1}\lambda^r c_\lambda^{e(G)}\zeta_\varphi^r ,\label{uppb-2}
    \end{eqnarray}
    where $x_{r+i}:=y_{i}$ for $i=1,\dots,m$. What's more, we can also bound $\kappa_2(N_G)$ from below
    \begin{eqnarray}
      \kappa_2(N_G)&\ge&\lambda^r c_{\lambda}^{e(G)} \int_{(\R^d)^r }
      \ \prod_{
        \substack{1 \leq k \textcolor{black}{<} l \le r+m
          \\
          \{ k , l \}\in E_G 
      }}\varphi(x_k,x_l) \ \! \mathrm{d}x_1 \cdots\mathrm{d}x_r\nonumber\\
      &=:&\lambda^r c_{\lambda}^{e(G)}C_4.
      \label{lowb-2}
    \end{eqnarray}
    Combining \eqref{uppb-2} and \eqref{lowb-2}, we get 
    \begin{align}
      \kappa_n(\widebar{N}_G)&\leq \frac{n!^rr!^{n-1}\lambda^r c_\lambda^{e(G)}\zeta_\varphi^r }{\big(\lambda^r c_{\lambda}^{e(G)}C_4\big)^{n/2}}\nonumber\\
      &\leq n!^r\big( \lambda^rc_\lambda^{e(G)}\big)^{-(n-2)/2}\frac{r!^{n-1}\zeta_\varphi^{rn}}{C_4^{n/2}}\nonumber\\
      &\leq \frac{n!^r}{\Big(C_5\sqrt{\lambda^r c_\lambda^{e(G)}}\Big)^{n-2}}
      \ ,            \end{align}
            where $C_5$ is a constant depending only on $r$ and $\varphi$.
  \end{enumerate}
  When $m=0$ the above arguments apply
  by replacing \textcolor{black}{the} upper bound \eqref{fjkldf2}
  with $\mu (\real^d)^{1+(r-1)n}$. 
\end{Proof}
\textcolor{black}{
We are now ready to prove Theorem~\ref{poslim}.
  \begin{Proofy} \hskip-0.16cm {\em of Theorem~\ref{poslim}}. 
  Since $\lim_{\lambda \to \infty} \lambda c_\lambda^{e(G)/r} = c >0$
  and $a_m(G) r \leq e(G)$,
  we have $ c_\lambda \lesssim \lambda^{-( r-1 )/(e(G)-a_m(G)) }$. 
    Hence, as in the proof of Corollary~\ref{themo-1-1}-\eqref{partc}, 
   by \eqref{cum-eq} we obtain 
    \begin{eqnarray*}
      \lim_{\lambda\to\infty}\kappa_n(N_G) & = &
      |\mathrm{Aut}_{\bullet}(G)|^{n-1}
      \lim_{\lambda\to\infty}\lambda^r 
      \int_{(\R^d)^r}
      \Bigg(
      \prod_{\textcolor{black}{\{i,j\}}\in E_G}
      c_\lambda \varphi(x_i,x_j)
      \Bigg)
      \ \! \mu(\mathrm{d}x_1)\cdots\mu(\mathrm{d}x_r)
      \\
      &=&|\mathrm{Aut}_{\bullet}(G)|^{n-1} c^r\int_{(\R^d)^r}
      \Bigg(
      \prod_{\textcolor{black}{\{i,j\}}\in E_G}\varphi(x_i,x_j)
      \Bigg) \ \!
      \mu(\mathrm{d}x_1)\cdots\mu(\mathrm{d}x_r), \qquad n \geq 1, 
    \end{eqnarray*}
  since the count of connected non-flat set partitions $\rho\in {\rm CNF}([n]\times[r])$ such that $|\rho|=r$ and $e(\rho_G)=e(G)$
  is $|\mathrm{Aut}_{\bullet}(G)|^{n-1}$. 
  We conclude from Theorem~6.14 in \cite{JLR}. 
\end{Proofy}
}

  \appendix

\section{Convex hull code}
\label{fjkldsf-11}
\noindent  
The following SageMath code determines the convex hull of $\Sigma_n(G,m)$ 
and its upper boundary $\widehat{\Sigma}_n(G,m)$,
see Figures~\ref{fig1}-\ref{fig2} and \ref{fig4}.
 This code and the following one are available for download at 
 \url{https://github.com/nprivaul/convex-hull}.  

\medskip
\smallskip

\begin{lstlisting}[language=Sage] 
def partitions(points):
      if len(points) == 1:
         yield [ points ]
         return
         first = points[0]
         for smaller in partitions(points[1:]):
         	for m, subset in enumerate(smaller):
			yield smaller[:m] + [[ first ] + subset]  + smaller[m+1:]
	yield [ [ first ] ] + smaller

def nonflat(partition,r):
  p = []
  for j in partition:    
      seq = list(map(lambda x: (x-1)//r,j))
      p.append(len(seq) == len(set(seq)))
  return all(p)

def connected(partition,n,r):
  q = []; c = 0
  if n  == 1: return all([len(j)==1 for j in partition])
  for j in partition:
      jk = list(set(map(lambda x: (x-1)//r,j)))
      if(len(jk)>1):            
          if c == 0:
              q = jk; c += 1
          elif(set(q) & set(jk)):
              d=[y for y in (q+jk) if y not in q]
              q = q + d
  return n == len(set(q))

def connectednonflat(n,r):
  points = list(range(1,n*r+1))
  randd = []
  for m, p in enumerate(partitions(points), 1):
      randd.append(sorted(p))
  for rou in range(r,(r-1)*n+2):    
      rs = [d for d in randd if (nonflat(d,r) and len(d)==rou)]
      rss = [e for e in rs if connected(e,n,r)]
      print("Connected non-flat partitions with",rou,"blocks:",len(rss))
  cnfp = [e for e in randd if (connected(e,n,r) and nonflat(e,r))]
  print("Connected non-flat set partitions:",len(cnfp))
  return cnfp

def graphs(G,EP,setpartition,n):
  r=len(set(flatten(G)));rhoG = []
  for j in range(n):
      for hop in G: rhoG.append([r*j+hop[0],r*j+hop[1]])
      for l in range(len(EP)):
          F=EP[l]
          for i in F: rhoG.append([j*r+i,n*r+l+1]);
  for i in setpartition:
      if(len(i)>1):
          b = []
          for j in rhoG:
              b.append([i[0] if ele in i else ele for ele in j])
          rhoG = b
  for i in rhoG: i.sort()
  return rhoG

def convexhull(n,G,EP):
  r=len(set(flatten(G)));m=len(EP)
  cnfp=connectednonflat(n,r)
  L=[]
  le=sum(len(EP[j]) for j in range(len(EP)))
  for setpartition in cnfp: 
      rhoG=graphs(G,EP,setpartition,n)
      edgesrhoG = [i for n, i in enumerate(rhoG) if i not in rhoG[:n]]
      vertrhoG = set(flatten(edgesrhoG));
      L.append((n*r-(len(vertrhoG)-m),n*(len(G)+le)-len(edgesrhoG)))
  return sorted(set(L))
\end{lstlisting}

\section{Graph counting code}
\label{fjkldsf-22}
\noindent 
The following R code uses the graph6 and sparse6 formats for undirected graphs
and the data files available at \url{https://users.cecs.anu.edu.au/~bdm/data/graphs.html}. 

\medskip
 
\begin{lstlisting}[language=R] 
library(rgraph6); library(igraph); library(matrixStats)
r=3; m=2; graphs=read_file6("graph5c.g6") # m+r=5
# r=5; m=3; graphs=read_file6("graph8c.g6") # m+r=8
count=0; total=0; trees=0; treesam=0; nontreesr=0; 
for (mat in graphs) {g=as.undirected(graph_from_adjacency_matrix(mat))
endpoints=combn(1:(m+r),m); lst = c();
for (k in 1:choose(m+r,m)) {
V(g)$color <- c(7, 2)[1 + V(g) %in% endpoints[,k]]
complement=setdiff(c(1:(m+r)),endpoints[,k])
if(sum(degree(subgraph(g,endpoints[,k])))==0 && is.connected(subgraph(g, complement)))
{if (all(sapply(lst, function(gg) !graph.isomorphic.vf2(g,gg)$iso))) {lst=c(lst,list(g)); 
count=count+1; a=0; exit=0; 
for (ii in complement) {a=max(a,sum(g[ii,endpoints[,k]]))}
for (i in (m+2):(m+r)){ for (j in 1:choose(m+r,i)){ 
h=subgraph(g, combn(1:(m+r),i)[,j])
if (((ecount(g)-a)/(vcount(g)-m-1))<((ecount(h)-a)/(vcount(h)-m-1))) {
exit=1;break;}}
if (exit==1) {break;}}
if (exit==0) {print(g); cat("a =",a,"\n"); total=total+1
trees=trees+is_tree(g)*(a==m); treesam=treesam+is_tree(g)
nontreesr=nontreesr+!is_tree(subgraph(g, complement))
plot(g); print("Working ...");}}}}}
cat("Tree with a<>m count = ",treesam, ";Tree count = ",trees,"out of total =",total,"out of", count, "\n");
\end{lstlisting}

\footnotesize

\newcommand{\etalchar}[1]{$^{#1}$}
\def\cprime{$'$} \def\polhk#1{\setbox0=\hbox{#1}{\ooalign{\hidewidth
  \lower1.5ex\hbox{`}\hidewidth\crcr\unhbox0}}}
  \def\polhk#1{\setbox0=\hbox{#1}{\ooalign{\hidewidth
  \lower1.5ex\hbox{`}\hidewidth\crcr\unhbox0}}} \def\cprime{$'$}

\end{document}